\documentclass[12pt,oneside]{amsart}
\title[An inverse problem for the Standard Model ]{An inverse problem for the Standard Model of particle physics}

\author[X. Chen]{Xi Chen}
\address{Shanghai Center for Mathematical Sciences, Fudan University, Shanghai 200438, China;	Shanghai Academy of Artificial Intelligence for Science, Shanghai 200232, China;
	Center for Applied Mathematics, Fudan University, Shanghai 200433, China. }
\email{xi\_chen@fudan.edu.cn}

\author[M. Lassas]{Matti Lassas}
\address{Department of Mathematics and Statistics, University of Helsinki, Helsinki FI-00014, Finland.  }
\email{Matti.Lassas@helsinki.fi}

\author[L. Oksanen]{Lauri Oksanen}
\address{Department of Mathematics and Statistics, University of Helsinki, Helsinki FI-00014, Finland.  }
\email{Lauri.oksanen@helsinki.fi}

\author[G.P. Paternain]{Gabriel P. Paternain}

\address{ Department of Mathematics, University
	of Washington, Seattle, WA 98195, USA.}
\email{gpp24@uw.edu}

\usepackage[asymmetric, centering]{geometry}
\usepackage{graphicx} 
\usepackage{mathrsfs}
\usepackage{amsfonts}
\usepackage{amscd}
\usepackage{amsmath}
\usepackage{hyperref}

\usepackage{enumitem}

\usepackage{tikz-cd,amsmath,amsfonts}

\IfFileExists{notused_foobar.sty}{\usepackage{tweakslo}}{
\usepackage{amssymb,thmtools,mathtools,todonotes}
\declaretheorem{theorem,definition,lemma,proposition,corollary,remark}}

\def\p{\partial}
\def\R{\mathbb R}
\DeclareMathOperator{\supp}{supp}

\usepackage{slashed} 

\def \M{{\mathbb M}}
\def\C{\mathbb C}
\def \g{{\mathfrak g}}
\def \J{{\mathbb J}}
\DeclareMathOperator{\Ad}{Ad}
\DeclareMathOperator{\vol}{vol}
\DeclareMathOperator{\id}{id}

\DeclareMathOperator{\GL}{GL}
\def\P{\mathbf P}
\def\V{\mathcal{V}}
\def\W{\mathcal{W}}

\allowdisplaybreaks[4]
\begin{document}
\begin{abstract} We pose and solve an inverse problem for the classical field equations that arise in the Standard Model of particle physics. Our main result describes natural conditions on the representations, so that it is possible to recover all the fields from measurements in a small set within a causal domain in Minkowski space. These conditions are satisfied for the representations arising in the Standard Model.
\end{abstract}
\maketitle

\tableofcontents

\section{Introduction} 
This paper discusses an inverse problem for the classical field equations of the Standard Model of particle physics where the objective is to recover the fields by performing measurements in a small set within a causal domain in Minkowski space, allowing waves to propagate and return.

The main geometric objects involved are spinor fields (sections of a twisted chiral spinor bundle), gauge fields (connections), and Higgs fields (sections of a suitable Hermitian vector bundle). Our focus here is on the spinor fields describing the matter content in the Standard Model, thus completing a long-term project that began with the work in \cite{CLOP1, CLOP2, CLOP3}. An overview of the sequence is as follows:

\begin{itemize}
	\item \cite{CLOP1} discusses a cubic wave equation for the connection wave operator, setting up the scene for our microlocal approach and proving a result on the broken (non-abelian) light ray transform that underpins the entire scheme.
	
	\item \cite{CLOP2} addresses the case of the pure Yang-Mills equations and the recovery of the gauge field.
	\item \cite{CLOP3} explores the Yang-Mills-Higgs equations and how to retrieve the Higgs field. An important insight here is to leverage the center $Z(\mathfrak{g})$ of the Lie algebra of the structure group $G$.
	\item The current paper introduces spinor fields, thus completing the full Lagrangian of the Standard Model, including Yukawa couplings that describe interactions between the Higgs field and spinors.
	
\end{itemize}

The particle content of the Standard Model is described by two representations $\rho$ (Higgs field) and $\varrho$ (fermions) of $G$ in Hermitian vector spaces. Our main result pinpoints precisely the conditions on $\rho$ and $\varrho$ that are needed for our microlocal approach to work, and we confirm that these conditions hold in the Standard Model.
Once again, the role of the center of the Lie algebra is crucial, as the approach works best when activating sources in $Z(\g)$. For the Standard Model, where 
	\begin{align*}
G = \mathrm{SU}(3) \times \mathrm{SU}(2) \times \mathrm{U}(1),
	\end{align*}
this corresponds to using the $\mathrm{U}(1)$ component and the fact that the hypercharge is non-zero for all fermions. Our main theorem is stated in detail as Theorem \ref{thm:main thm} below for an arbitrary compact matrix Lie group and general representations $\rho$ and $\varrho$.

A direct approach to the main theorem activating all sources and expecting that the nonlinear interactions will deliver the recovery is infeasible, as the full system of equations is too complex. The goal of our quartet of papers is to explain how to recover each field in an orderly fashion, starting with the connection, then the Higgs field, and finally the spinors. Nonetheless, the computations required are substantial; in \cite{CLOP2} certain parts were checked by a computer code, while \cite{CLOP3} was approachable by hand but still challenging. In the current paper, we also use computer assistance for the final computations involving the gamma matrices of the Clifford structure. Elucidating the role of the latter structure in the recovery of the spinor fields is one of the main aims of the present paper.

\subsection{The Lagrangian}

\subsubsection{Minkowski space} 

Throughout the paper we work in Minkowski space $(\M, g)$, the usual background for the Standard Model. We fix the signature $(-,+,+,+)$ and adopt Cartesian coordinates $(x^0, x^1, x^2, x^3)$ on $\M$ such that $$g = - dx^0 \otimes dx^0 + dx^1 \otimes dx^1 + dx^2 \otimes dx^2  + dx^3 \otimes dx^3.  $$


\subsubsection{Clifford algebra, gamma matrices and Dirac form} We shall  work with the Clifford algebra of Minkowski space exclusively via the physical gamma matrices that
determine the chiral (or Weyl) representation. They are given by 
\begin{align*}\Gamma_0 &:= \imath\left(  \begin{array}{cc}
		0 & \id_2 \\
		\id_2 & 0
	\end{array}  \right),   \\  
	\Gamma_k &:= \imath\left(  \begin{array}{cc}
		0 & \sigma_k \\
		-\sigma_k & 0
	\end{array}  \right) \quad k = 1, 2, 3,
\end{align*} where  $\imath = \sqrt{-1}$ and $\sigma_k$ are the Pauli matrices 
\[\sigma_1 = \left(\begin{array}{cc}  0 & 1 \\ 1 & 0 \end{array}\right), \quad \sigma_2 = \left(\begin{array}{cc}  0 & - \imath \\ \imath & 0 \end{array} \right), \quad \sigma_3 = \left(\begin{array}{cc}  1 & 0 \\ 0 & -1 \end{array} \right).\]
We will raise and lower indices using the Minkowski metric. In particular, 
	\begin{align*}
\Gamma^\alpha=g^{\alpha \beta}\Gamma_{\beta}
	\end{align*}
and the key identity giving rise to the Clifford structure is
	\begin{align}\label{commrel_Gamma}
\Gamma^{\alpha}\Gamma^{\beta}+\Gamma^\beta \Gamma^\alpha=2g^{\alpha\beta}\id_{4}.
	\end{align}
(Switching the signature to $(+,-,-,-)$ changes the gamma matrices as described in \cite[Example 6.3.18]{Ha}.)
The gamma matrices give rise to a Clifford multiplication on $\Delta:=\C^{4}$:
\begin{align*}\M\times \Delta&\longrightarrow \Delta\\(X,\psi)&\longmapsto X\cdot \psi:=\imath X^{a}\Gamma_{a}\psi\end{align*}
where $X=X^{a}e_{a}$ and $e_{a}$ is a canonical basis vector in $\M$. There is also a chirality operator given by
\[\Gamma_{5}:=-\imath \Gamma_0 \Gamma_{1}\Gamma_{2} \Gamma_3,\]
that decomposes $\Delta$ into its $\pm 1$-eigenspaces denoted by $\Delta_{L}$ and $\Delta_{R}$. Below, we will often write $\psi\in\Delta$ as
$(\psi_{L},\psi_{R})$ using this splitting (which simply separates the first two coordinates in $\C^4$ from the last two).

We will endow $\Delta$ with a specific non-degenerate $\R$-bilinear form (called a Dirac form)
\[\langle \cdot, \cdot \rangle_{\Delta}: \Delta\times \Delta \longrightarrow \C\]
given by 
\[\langle \psi, \phi \rangle_{\Delta} =\psi^{\dagger} (-\imath\Gamma_{0})\phi=\psi^{\dagger}_{L}\phi_{R}+\psi_{R}^{\dagger}\phi_{L}\]
where we think of $\psi,\phi\in\Delta$ as column vectors in $\C^4$ and $\dagger$ denotes conjugate transpose.
It has the properties:
\begin{enumerate}[label={\arabic*.}]
	\item $\langle X \cdot \psi, \phi \rangle_{\Delta} = - \langle \psi, X \cdot \phi \rangle_{\Delta}$ for all $ X \in \M$ and $ \psi, \phi \in \Delta$;
	\item $\langle \psi, \phi \rangle_{\Delta} = \overline{\langle \phi, \psi \rangle_{\Delta}}$ for all $ \psi, \phi \in {\Delta}$;
	\item $\langle \psi, c\phi \rangle_{\Delta} = c \langle \psi, \phi \rangle_{\Delta} = \langle \bar{c} \psi, \phi \rangle_{\Delta}$ for all  $\psi, \phi \in {\Delta}$ and $c \in \C$.
	
\end{enumerate}

Note that the bilinear form is null on $\Delta_{L}$ and $\Delta_{R}$ and
\[\langle \psi, \phi \rangle_{\Delta} =\langle \psi_{L}, \phi_{R} \rangle_{\Delta} +\langle \psi_{R}, \phi_{L} \rangle_{\Delta}.\]



\subsubsection{Representations} Let $G$  be a compact, connected matrix Lie group and endow its Lie algebra $\g$ with an $\text{Ad}$-invariant inner product.
We consider $\V$ and $\W$ two finite dimensional complex vector spaces, equipped with Hermitian inner products $\langle \cdot, \cdot \rangle_{\V}$ and $\langle \cdot, \cdot \rangle_{\W}$. Suppose we are given two complex linear representations
\begin{align*}
	\varrho:G\to \text{GL}(\V), \quad	\rho:G\to \text{GL}(\W),
\end{align*}
leaving invariant the respective Hermitian inner products.  We shall moreover, assume an orthogonal splitting  
\begin{align*} 
	\V:=\V_{L}\oplus \V_{R} 
\end{align*}
as well as a splitting of the representation $\varrho=\varrho_{L}\oplus \varrho_{R}$.  This assumption is to account for the chirality of the matter content in the Standard Model.

\subsubsection{Bundles and sections of interest} The fields of interest will be sections of suitable vector bundles. Since we are working on Minkowski space, we shall only consider trivial bundles. They are:

\begin{align*} 
	S := \M \times \Delta, \quad F := \M\times \V,
	\quad\Ad := \M\times \g,\quad 
	E := \M\times \W.
\end{align*} The bundle metrics for these bundles  are induced from the Dirac form $\langle \cdot, \cdot \rangle_{\Delta}$, the Hermitian inner products on the fibres of $F$ and $E$, and the $\text{Ad}$-invariant inner product on the Lie algebra $\g$ respectively.

The spinor bundle $S$ splits into a direct sum of complex Weyl spinor bundles $$S = S_L \oplus S_R,\quad S_L := \M \times \Delta_L \quad S_R := \M \times \Delta_R.$$

Next, we introduce the twisted spinor bundle $$S\otimes F = \M \times \left(\Delta \otimes \V\right)$$ and the twisted chiral spinor bundles
\begin{align*}\left(S \otimes F\right)_+ &:= \left(S_L \otimes F_L\right) \oplus \left(S_R \otimes F_R\right) = \M \times  (\Delta \otimes \V)_+   \end{align*} where
\begin{align*}
	(\Delta \otimes \V)_+ &:= \left(\Delta_L \otimes \V_L\right) \oplus \left(\Delta_R \otimes \V_R\right).
\end{align*}
and 
\begin{align*}\left(S \otimes F\right)_{-} &:= \left(S_R \otimes F_L\right) \oplus \left(S_L \otimes F_R\right) = \M \times  (\Delta \otimes \V)_-  \end{align*} where
\begin{align*}
	(\Delta \otimes \V)_- &:= \left(\Delta_R \otimes \V_L\right) \oplus \left(\Delta_L \otimes \V_R\right).
\end{align*}

We will be working with spaces of functions
$$ C^\infty(\M, \Delta \otimes \V), \quad C^\infty(\M, \g),\quad C^\infty(\M, \W),$$ and the spaces of $\g$-valued $k$-forms 
\begin{align*} 
	\Omega^m(\M, \Delta \otimes \V) &:= C^\infty(\M, \Delta \otimes \V) \otimes \Omega^m(\M, \C), \\ 
	\Omega^m(\M, \g) &:= C^\infty(\M, \g) \otimes \Omega^m(\M, \C),\\  
	\Omega^m(\M, \W) &:= C^\infty(\M, \W) \otimes \Omega^m(\M, \C).\end{align*}
The inner product of each vector space along with the wedge product of differential forms gives a natural pairing for vector-valued forms. Specifically, we have
\begin{align*}
	\wedge_{\Delta \otimes \V} : \Omega^p(\M, \Delta\otimes\V) \times \Omega^q(\M, \Delta\otimes\V) &\longrightarrow \Omega^{p + q}(\M, \C),\\
	\wedge_{\g} : \Omega^p(\M, \g) \times \Omega^q(\M, \g) &\longrightarrow \Omega^{p + q}(\M, \R), \\
	\wedge_{\W} : \Omega^p(\M, \W) \times \Omega^q(\M, \W) &\longrightarrow \Omega^{p + q}(\M, \C).
\end{align*} The $L^2$-inner products of such forms are defined by 
\begin{align*}		\langle \alpha, \beta \rangle_{L^2, S \otimes F} &:= \int_\M \alpha \wedge_{\Delta \otimes \V} \star \beta  && \forall \alpha, \beta \in \Omega^m(\M, (\Delta \otimes \V)_+),\\
	\langle \alpha, \beta \rangle_{L^2, \Ad} &:= \int_\M \alpha \wedge_\g \star \beta  && \forall \alpha, \beta \in \Omega^m(\M, \g),\\ 
	\langle \alpha, \beta \rangle_{L^2, E} &:= \int_\M \alpha \wedge_\W \star \beta  && \forall \alpha, \beta \in \Omega^m(\M, \W),
\end{align*} where $\star$ denotes the Hodge star operator on $\M$. 

The Standard Model describes fermions with twisted chiral spinor fields $\psi \in C^{\infty}(\M, (\Delta \otimes \V)_+)$, 
gauge bosons with Yang--Mills gauge fields $A \in \Omega^{1}(\M,\g)$, and bosons via the Higgs fields $\Phi \in C^{\infty}(\M,\mathcal{W})$.

\subsubsection{Differential operators}\label{sec_diff_ops}  Given $A\in \Omega^{1}(\M,\g)$, the induced exterior covariant derivatives in $\Ad$ and $E$ are defined  respectively to be
\begin{eqnarray*}
	D_A : \Omega^m(\M, \g) &\longrightarrow& \Omega^{k+1}(\M, \g) \\  W &\longmapsto& dW + [A,   W],
	\\
	d_A : \Omega^m(\M, \mathcal W) &\longrightarrow& \Omega^{k+1}(\M, \mathcal W) \\  \Upsilon &\longmapsto& d\Upsilon + \rho_\ast(A) \wedge \Upsilon.
\end{eqnarray*}
The covariant derivatives $D_A$ and $d_{A}$ are compatible with the inner products in $\g$ and $\mathcal{W}$ respectively.
The formal adjoints of $D_A$ and $d_A$ are defined through the Hodge star operator $\star$   to be 
\begin{align*}
	D_A^\ast =  \star D_A \star, \quad	d_A^\ast = \star d_A \star.
\end{align*}

Analogously, the twisted spin covariant derivative of a smooth section $\psi$ on the twisted bundle $(S \otimes F)_+$ takes the form 
\begin{align*}d_A \psi := d \psi + \varrho_\ast (A) \psi.\end{align*} 
The twisted Dirac operator on $(S \otimes F)_+$ is given by 
\begin{align*}\slashed{D}_A : C^\infty(\M,(\Delta\otimes \V)_{+}) & \longrightarrow C^\infty(\M,(\Delta\otimes \V)_{-})\\  \psi &\longmapsto \imath \Gamma^a \left(d\psi(e_a) + \varrho_\ast (A_a) \psi\right)=\imath \Gamma^{a}\,d_{A}\psi(e_{a}).\end{align*} 
If we decompose the spinor field as \[\psi : = \psi_L + \psi_R \in C^\infty(\M, \Delta_L \otimes \V_L) + C^\infty(\M, \Delta_R \otimes \V_R),\]  the operator splits accordingly into
\begin{align*}\slashed{D}_A \psi &:= \slashed{D}_A \psi_L + \slashed{D}_A \psi_R,  \\  
	\slashed{D}_A \psi_{L}	&:= \imath \Gamma^a \left(d\psi_{L}(e_a) + \varrho_{L, \ast} (A_a) \psi_{L}\right),  \\  
	\slashed{D}_A \psi_{R}	&:= \imath \Gamma^a \left(d\psi_{R}(e_a) + \varrho_{R, \ast} (A_a) \psi_{R}\right).
\end{align*}
The twisted Dirac operator is formally self-adjoint, i.e. 
\[\langle \slashed{D}_A \phi, \psi \rangle_{S \otimes F, L^2} = \langle  \phi, \slashed{D}_A \psi \rangle_{S \otimes F, L^2}.\]
We also recall that the curvature form $F_A$ of $A$ is a $\g$-valued 2-form
\begin{align}
    \label{curvature form FA}
    F_A=dA+\frac 12 [A,A].
\end{align}

\subsubsection{Yukawa coupling} 
Fermions acquire mass through interactions with Higgs bosons. This is incorporated in the Lagrangian via the Yukawa coupling (see \cite[Section 7.7]{Ha}).

The Yukawa form is a map
\[\mathrm{Y}:\V_{L}\times \W\times \V_{R}\to\C\]
which is invariant under the action of $G$, complex antilinear in $\V_{L}$, real linear in $\W$ and complex linear in $\V_{R}$.

Given a Yukawa form $\mathrm{Y}$ and a constant $g_{\mathrm{Y}}$   we define the Yukawa coupling 
\begin{align*}\mathbf{Y}:(\Delta_L \otimes \V_{L})\times \W\times (\Delta_R\otimes \V_{R}) &\longrightarrow \mathbb{R}\\ 
	(s_L\otimes \tau_{L},\phi, s_R\otimes \tau_{R}) & \longmapsto -2g_{\mathrm{Y}}\Re(\langle s_L, s_R\rangle_{S}\,\mathrm{Y}(\tau_{L},\phi,\tau_{R})).\end{align*}

\subsubsection{The Lagrangian} 
We have now all the necessary ingredients to write down the Lagrangian density (4-form) of the Standard Model:
\begin{align*} {\mathscr{L}_{\mathrm{SM}}(\psi_L, \psi_R, A, \Phi) : =    \mathscr{L}_{\mathrm{D}}[A, \psi_L, \psi_R] + \mathscr{L}_{\mathrm{Y}}[\psi_L, \Phi, \psi_R]} + \mathscr{L}_{\mathrm{YM}}[A] + \mathscr{L}_{\mathrm{H}}[\Phi, A]\end{align*} 
where 
\begin{align*} 
	\mathscr{L}_{\mathrm{D}}[\psi_L, \psi_R, A] &:= (\Re  \langle \psi_L, \slashed{D}_A \psi_L \rangle_{\Delta \otimes \V} +  \Re \langle    \psi_R, \slashed{D}_A \psi_R  \rangle_{\Delta \otimes \V} )\,d\vol   \\ 
	\mathscr{L}_{\mathrm{Y}}[\psi_L, \psi_R, \Phi] &:= \mathbf{Y}\left(\psi_L, \Phi, \psi_R\right)\,d\vol
	\\
	\mathscr{L}_{\mathrm{YM}}[A] &:= - \frac{1}{2}  F_A \wedge_{\g} \star F_A   \\ 
	\mathscr{L}_{\mathrm{H}}[A, \Phi] &:=    d_A \Phi \wedge_{\W} \star d_A \Phi  -  \mathbf{V}(|\Phi|_{\W}^2) \,d\vol  .\end{align*}
The Higgs potential $\mathbf{V}(|\Phi|_{\W}^2)$ is a quartic polynomial in $\Phi$. However, the specific coefficients of $\mathbf{V}$ do not play a role in the proofs. For simplicity, we take  \[\mathbf{V}(|\Phi|^2_{\W}) =  \frac{1}{2} |\Phi|_{\W}^4  - |\Phi|_{\W}^2.\]

The action functional of the Standard Model is the integration of the lagrangian over $\M$
\begin{align*}  \mathscr{A}_{\mathrm{SM}}(\psi, A, \Phi) : =  \int_{\M} \mathscr{L}_{\mathrm{SM}}(\psi, A, \Phi)  .\end{align*} We write $\psi \in C^\infty(\M, (\Delta \otimes \V)_+)$ instead of $\psi_L + \psi_R$ when there is no need to split the spinor field.

\subsection{The inverse problem}
We are interested in the critical points of the action functional $\mathscr{A}_{\mathrm{SM}}$ which are triples $(\psi,A,\Phi)$ lying in
\[ C^\infty(\M, (\Delta \otimes \V)_+) \times \Omega^1(\M, \g) \times C^\infty(\M, \W).\]  
The action functional has a well known gauge invariance and hence the set of critical points remains invariant under the action of the gauge group. Our result needs to account for this invariance.


\begin{figure}
	\centering
	\includegraphics[width=0.5\textwidth,trim={3cm 5cm 6cm 2cm},clip]{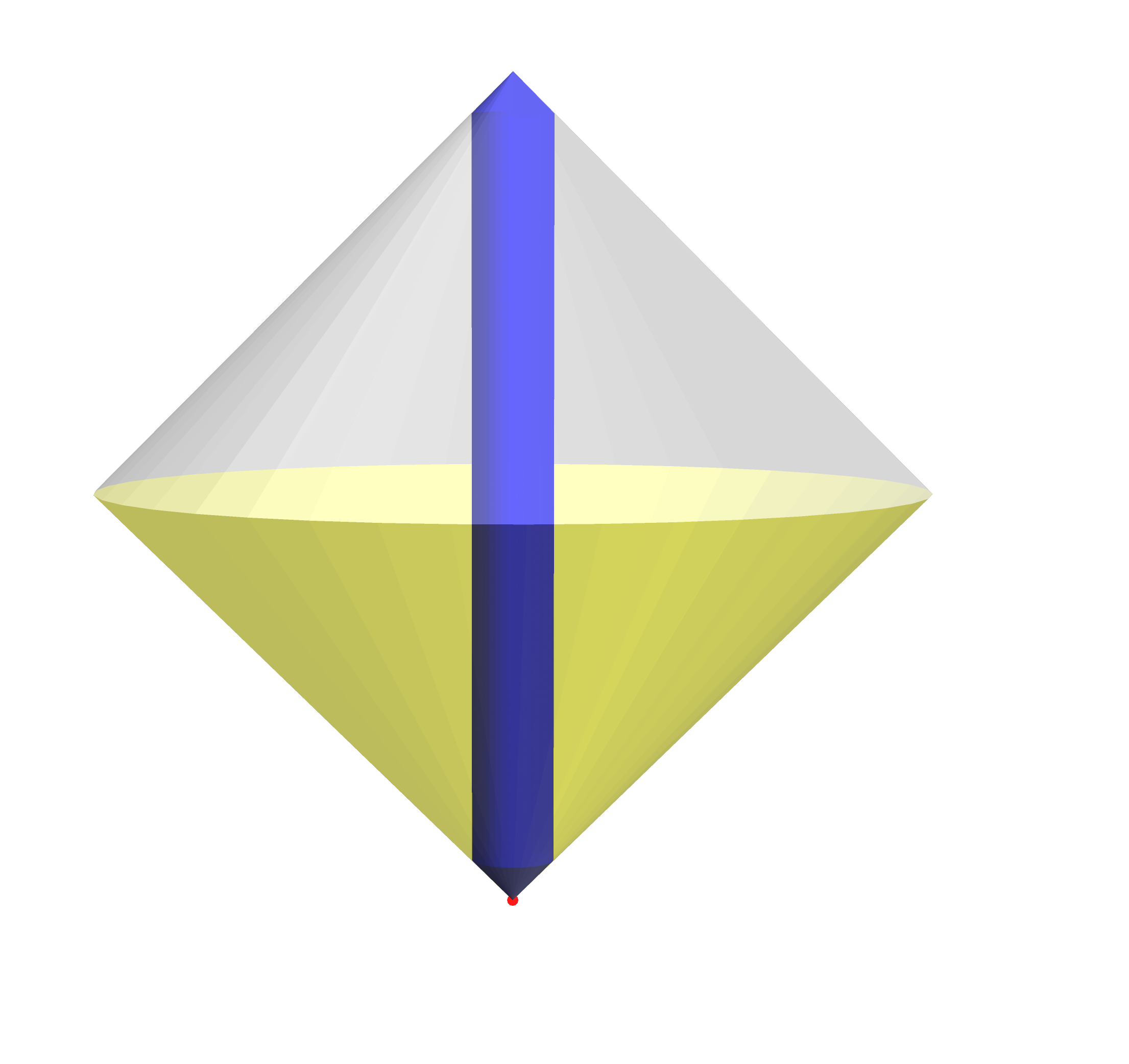}
	\caption{The set $\mho$ (in blue) inside the diamond $\mathbb D$ in the $1+2$ dimensional case. The part $\p^- \mathbb D$ of the boundary of $\mathbb D$ is shaded in yellow. The point $p$ is drawn in red.}
	\label{fig_D}
\end{figure}

The set up involves the causal diamond (see Figure \ref{fig_D}):
	\begin{align}\label{def_diamond}
\mathbb D
= \{ (t,x) \in \M : |x| \le t + 1,\ |x| \le 1 - t \}.
	\end{align}
We will show that a critical triple $(\psi, A, \Phi)$ can be recovered from suitable data in the domain
\begin{align}\label{def_mho}
	\mho = \{(t,x) : \text{$(t,x)$ is in the interior of $\mathbb D$ and $|x| < \varepsilon_0$} \},
\end{align}
($0 < \epsilon_0 < 1$) modulo the action of the gauge group  $C^\infty (\mathbb{D}; G)$. For technical reasons, it is preferable
to consider the pointed gauge group \[G^0(\mathbb{D}, p) = \{\mathbf{U} \in C^\infty (\mathbb{D}; G) : \mathbf{U}(p) = \id\}, \quad \mbox{for $p = (-1, 0) \in \mathbb{D}$}.\] 
An element $\mathbf{U}$ of $G^0(\mathbb{D}; G)$ acts on triples $(\psi, A, \Phi)$ from the right as follows: 
\begin{align}\label{eqn : gauge transform}(\psi, A, \Phi) \cdot \mathbf{U} = \left(s \otimes \left( \varrho(\mathbf{U}^{-1})\tau\right), A \cdot \mathbf{U}, \rho(\mathbf{U}^{-1}) \Phi \right)\end{align} 
where we write the twisted spinor field in the form $\psi := s \otimes \tau$ and then extend by linearity.

We write $\sim$ for the equivalence of two triples under the action of an element in $G^0(\mathbb{D}, p)$; in other words:
\begin{align*} (\psi, A, \Phi) \sim (\zeta, B, \Xi) \Longleftrightarrow 
	\exists\, \mathbf{U} \in G^0(\mathbb{D}, p)\; \mbox{such that}\, (\psi, A, \Phi) \cdot \mathbf{U} = (\zeta, B, \Xi)\end{align*}

We take the data set of a critical point $(\psi, A, \Phi)$ in $\mho$ to be
\begin{align*}
    \mathcal{D}_{(\psi, A, \Phi)}(\mho): = \left\{\begin{array}{l|l}
	(\phi, V, \Psi)|_\mho & 
	\left.\begin{array}{l}
		\mbox{$(\phi, V, \Psi) \in C^3(\mathbb{D})$ is a critical point of $\mathscr{A}_{\mathrm{SM}}$ } \\ \mbox{in $\mathbb{D} \setminus \mho$ and $(\phi, V, \Psi)\sim (\psi, A, \Phi)$ near $\partial^- \mathbb{D}$} \end{array} \right.  \end{array} \right\}.
        \end{align*}
Here
\[	\partial^- \mathbb{D}=\{(t,x)\in\mathbb{D}:\,|x|=t+1\}\]
and 	
$(\phi, V, \Psi) \sim (\psi, A, \Phi)$ near $\partial^- \mathbb{D}$ if there are $\mathbf{U} \in G^0(\mathbb{D}, p)$ and a neighbourhood $\mathcal  U$ of $\partial^- \mathbb{D}$ such that $(\phi, V,\Psi) = (\psi, A,\Phi)\cdot \mathbf{U}$ on $ \mathcal {U}\cap \mathbb{D}$.
Note that $p = (-1, 0) \in \bar{\mho}$.

{Let us describe the data set $\mathcal{D}_{(\psi, A, \Phi)}(\mho)$ in alternative terms.
The property that $(\phi, V, \Psi) \in C^3(\mathbb{D})$ is a critical point of $\mathscr{A}_{\mathrm{SM}}$ in $\mathbb{D} \setminus \mho$ is equivalent to 
the fields $(\phi, V, \Psi)$ satisfying a system of 
non-linear partial differential
equations 
\begin{align}\label{standard short hand2}
    \mathcal P(\phi, V, \Psi
    )=(\mathcal{K}, \mathcal{J}, \mathcal{F})\quad\hbox{in }
    \mathbb{D}
\end{align}
with some sources $(\mathcal{K}, \mathcal{J}, \mathcal{F})$ that are supported
in the closure $\overline  \mho$. 
Here \eqref{standard short hand2} is shorthand for the system of equations \eqref{pL}-\eqref{pH} below.
When it holds,
the sources $(\mathcal{K}, \mathcal{J}, \mathcal{F})$ have to 
satisfy the compatibility condition \eqref{eqn : compatibility}. 
In addition, we require that $(\phi, V, \Psi)$ satisfies the
initial conditions 
\begin{align}\label{standard short initial}
(\phi, V, \Psi)\sim (\psi, A, \Phi),\quad \hbox{near $\partial^- \mathbb{D}$},
\end{align}
where $(\psi, A, \Phi)$ is an unperturbed, background solution of the Standard Model, that is, it solves \eqref{standard short hand2} with a vanishing right-hand side.
Then $\mathcal{D}_{(\psi, A, \Phi)}(\mho)$
can be rewritten as the following set of source-solution pairs for 
the initial 
        value problem \eqref{standard short hand2}-\eqref{standard short initial} 
\begin{align}\label{data set D2}
   \left\{\begin{array}{l|l}
	\big((\mathcal{K}, \mathcal{J}, \mathcal{F}),(\phi, V, \Psi)|_\mho\big) & 
	\left.\begin{array}{l}
		\mbox{$(\phi, V, \Psi) \in C^3(\mathbb{D})$ satisfy \eqref{standard short hand2}-\eqref{standard short initial},} \\ \mbox{$\supp(\mathcal{K}, \mathcal{J}, \mathcal{F})\subset\overline \mho$}
        \end{array} \right.  \end{array} \right\}.
        \end{align}
}
\medskip

Under suitable assumptions on the Lie algebra $\g$ and the derivatives at the identity of the representations $\rho$ and $\varrho$ we shall prove that one can uniquely determine the gauge equivalence class of  $(\psi, A, \Phi)$ in $\mathbb{D}$ from the local information $\mathcal{D}_{(\psi, A, \Phi)}(\mho)$ in $\mho$. Let $Z(\g)$ denote the centre of the Lie algebra $\g$.

\begin{definition} We shall say that a representation $\rho$ is \text{\bf hypercharged} if there exists $X\in Z(\g)$ such that
	\[\det \rho_{*}(X)\neq 0\]
	where $\rho_{*}:\g\to\text{\rm End}(\mathcal{W})$ denotes the derivative of $\rho$.
	\label{def:hyper}
\end{definition}

\begin{remark} {\rm Note that being hypercharged implies in particular that the centre $Z(\g)$ is non-trivial. As we will see below we are mostly interested in the case
		where $\dim Z(\g)=1$.}
\end{remark}


\begin{theorem}\label{thm:main thm}  Assume $\rho$ and $\varrho$ are hypercharged and that $Z(\g)\cap \text{\rm Ker}\,\rho_{*}=\{0\}$.
	Let $(\psi, A, \Phi)$ and $(\zeta, B, \Xi)$ be two critical points of $\mathscr{A}_{\mathrm{SM}}$ in $\mathbb{D}$.
	Then $$\mathcal{D}_{(\psi, A, \Phi)}(\mho) = \mathcal{D}_{(\zeta, B, \Xi)}(\mho) \Longleftrightarrow (\psi, A, \Phi)\sim (\zeta, B, \Xi) \, \mbox{in $\mathbb{D}$}.$$
\end{theorem}

\begin{remark} 
{\rm In our previous work \cite{CLOP3}
we proved the analogue of Theorem \ref{thm:main thm} for the Yang-Mills-Higgs system where there are no spinor fields. 
In \cite{CLOP3} we assumed a weaker condition on the representation $\rho$; however, the proof works as written only under the stronger hypercharged condition. 
Both conditions agree when $\dim Z(\g)=1$, so for the case of the Standard Model the distinction is immaterial.

}
\end{remark}

There is a rather general theory of inverse problems for nonlinear real principal type operators \cite{oksanen2024}. While wave type operators are included as a special case, this theory differs from the problem studied in this paper in two crucial ways. In the real principal type setting, the direct problem lacks gauge invariance, and the goal is to recover coefficients in the equation. In contrast, for the Standard Model equations, gauge invariance is present, and the unknowns to be determined are solutions of the equations rather than independent coefficients. Notably, inverse problems for the Einstein equations \cite{KLOU, Uhlmann2018} share these two features.

In fact, inverse problems for the equations of the Standard Model and Einstein equations are more akin to unique continuation problems than to coefficient determination problems. However, classical unique continuation results \cite{Hormander-Vol4} are inapplicable here, as they are known to fail in geometric settings lacking pseudoconvexity \cite{Alinhac1983}. Unique continuation from $\mho$ to $\mathbb{D}$ exemplifies such a setting. Alternative unique continuation results that do not require pseudoconvexity, such as \cite{Tataru}, impose analyticity conditions and are therefore unsuitable for our purposes.

\subsection{Representations for the Standard Model} The hypotheses on $\g$, $\rho$ and $\varrho$ are general enough to cover the case of the Standard Model, where 
\[G=\mathrm{SU}(3)_{C}\times \mathrm{SU}(2)_{L}\times \mathrm{U}(1)_{Y}\]
and the subindices $C,L$ and $Y$ stand for colour, weak and hypercharge respectively. The group $G$ acts on the Higgs vector space $\W=\C^2$ (equipped with the standard Hermitian inner product) as
\[\rho(g,h,e^{\imath\theta})w=e^{\imath n_{Y}\theta}hw,\]
where $(g,h, e^{\imath\theta})\in \mathrm{SU}(3)_{C}\times \mathrm{SU}(2)_{L}\times \mathrm{U}(1)_{Y}$ and the usual convention is that $n_{Y}=3$. We see that the factor $\mathrm{SU}(3)_{C}$ acts trivially on $\W$, 
$\text{Ker}\,\rho_{*}=\mathfrak{su}(3)$ and since $n_{Y}=3$ the representation is hypercharged ($\dim Z(\g)=1$). Thus the hypotheses on $\rho$ are satisfied.


The matter content of the Standard Model (fermions) is captured by the representation $\varrho$. This representation is highly reducible and breaks down into
several core components as we now briefly describe; details can be found in \cite{BH_10} and \cite[Chapter 8]{Ha}.
The first splitting is of the form $\V=\V_{L}\oplus \V_{R}$ with $\dim \V_{L}=24$ and $\dim \V_{R}=21$. Each $\V_{L}$ and $\V_{R}$ splits into three generations
and each generation splits into quark and lepton sectors. The core representations $Q^{i}_{L,R}$ and $L_{L,R}^{i}$ ($i=1,2,3$) 
appear in \cite[Table 8.1]{Ha} and they are actually quite simple; see also \cite[Table 1]{BH_10}. They are described as outer tensor products of each factor of $G$. The only representations that appear are the fundamental ones on $\C^{3}$ and $\C^2$ for $\mathrm{SU}(3)_{C}$ and  $\mathrm{SU}(2)_{L}$ respectively, the trivial representation, and for $\mathrm{U}(1)_{Y}$ we have the actions on $\C$ given by
	\begin{align}\label{def_Cy}
(e^{i\theta},z)\longmapsto e^{3yi\theta}z
	\end{align}
for some $y \in \R$. Moreover, \cite[Table 8.1]{Ha} shows that each core representation contains a factor of the form \eqref{def_Cy} with $y\neq 0$, and this is enough to guarantee that the resulting final representation $\varrho$ is hypercharged. It is because of the role played by $\mathrm{U}(1)_{Y}$ that we have chosen the name ``hypercharged'' in Definition \ref{def:hyper}.

The discussion above does not include antiparticles, but these can be incorporated by considering the dual representation on $\V^*$. Note that if we include the hypothetical right-handed neutrinos (or inert neutrinos), then this would violate the property of being hypercharged and our Theorem \ref{thm:main thm}  would not be able to detect them.
This is hardly surprising as these hypothetical particles do not interact with any of the gauge bosons of the Standard Model.

Finally we note that as byproduct of our results, the actual explicit form of the Yukawa coupling is irrelevant for the proof of Theorem \ref{thm:main thm}. While this term in the Lagrangian is crucial for the interactions between fermions and the Higgs boson, it is not relevant for our microlocal computation which is solidly anchored on the $\mathrm{U}(1)_{Y}$-factor.

\subsection{Structure of the proof}

Since the critical points of the action functional $\mathscr{A}_{\mathrm{SM}}$ satisfy the Euler-Lagrange equations \eqref{L}-\eqref{H}, we reformulate Theorem \ref{thm:main thm} as an inverse coefficient problem for a partial differential equation (PDE). The appropriate PDE for this purpose is the Euler-Lagrange equation with sources, written in terms of perturbed fields in the system of equations \eqref{eq:perturbed D+ prelim}-\eqref{eq:perturbed H prelim} below. The underlying fields to be recovered $(\psi, A, \Phi)$ appear as the coefficients of \eqref{eq:perturbed D+ prelim}-\eqref{eq:perturbed H prelim}.

Unlike conventional PDE inverse problems, several challenges arise in defining a source-to-solution map for the system \eqref{eq:perturbed D+ prelim}-\eqref{eq:perturbed H prelim}:
\begin{itemize}
\item The system is constrained by the compatibility condition \eqref{eqn : compatibility}, derived in Section \ref{subsection : compatibility}.
\item The Dirac equations \eqref{eq:perturbed D+ prelim}-\eqref{eq:perturbed D- prelim} differ from the wave type equations \eqref{eq:perturbed YM prelim}-\eqref{eq:perturbed H prelim},
calling for the conversion of the former into wave equations \eqref{eq:perturbed D+ convt}-\eqref{eq:perturbed D- convt}, as discussed in Section \ref{subsec : hyperbolic}.
\item Gauge conditions play a crucial role: the direct problem behaves well in the relative Lorenz gauge, and temporal gauge condition allows for relating the data set $\mathcal{D}_{(\psi, A, \Phi)}(\mho)$ to the direct problem. These two gauge conditions are considered in Sections \ref{subsec : Lorenz}-\ref{subsec : temporal}.
\item With these refinements, Section \ref{subsec : StS} establishes a source-to-solution map $\mathbf{L}_{(\psi, A, \Phi)}$ and reduces Theorem \ref{thm:main thm}  to Theorem \ref{thm : StS to fields}.
\end{itemize}

Having the source-to-solution map at hand, the remaining task is to solve the inverse coefficient problem for the system \eqref{eq:perturbed D+ prelim}-\eqref{eq:perturbed H prelim} using knowledge of $\mathbf{L}{(\psi, A, \Phi)}$:
\begin{itemize}
\item Section \ref{subsec : linear dyn} analyzes the dynamics of a linearized version 
\eqref{eq:SM pD+ linearized one-fold}-\eqref{eq:SM pH linearized one-fold} of the system.
\item Sections \ref{subsec : linearization}-\ref{subsec : from nonlinear to multilinear} use
these linear  propagation properties together with the technique of multiple linearizations. This technique originates in \cite{KLU}, where it was used to recover the leading-order terms in a semilinear wave equation with a simple power-type nonlinearity. Here it allows us 
to extract principal symbols of certain special solutions to the three-fold linearized system \eqref{eq:SM pD+ linearized 3-fold}-\eqref{eq:SM pH linearized 3-fold} from $\mathbf{L}{(\psi, A, \Phi)}$.
\end{itemize}

The final inversion is performed in two steps:
\begin{enumerate}[label={\arabic*.}]
\item At the principal level, the Yang--Mills and Higgs \eqref{eq:SM pYM linearized 3-fold}-\eqref{eq:SM pH linearized 3-fold} decouple from the Dirac equations \eqref{eq:SM pD+ linearized 3-fold}-\eqref{eq:SM pD- linearized 3-fold}. Therefore, Section \ref{subsec : YM recovery} first reconstructs the Yang--Mills and Higgs components by solving only the Yang--Mills and Higgs equations, employing the approach from \cite{CLOP3}.
\item For the spinor fields, Section \ref{subsec : EM signal in linearized system} introduces sources located exclusively in the center of the Lie algebra $Z(\g)$ and we exploit the Clifford algebra structure for the final recovery. This requires a delicate microlocal computation and makes use of the fact that $\varrho$ is hypercharged.

\end{enumerate}



\subsection{Structure of the paper} The remainder of the paper is organized as follows. Section \ref{sec:prelim} contains preliminaries, including the field equations, and the various bilinear forms that arise. Section \ref{sec:sts} explains how to transition from the data set $\mathcal{D}_{(\psi, A, \Phi)}(\mho)$ to a source-to-solution map $\mathbf{L}_{(\psi, A, \Phi)}$. Section \ref{sec : linearization} performs the linearizations and the final Section \ref{sec:recovery} recovers all the fields and proves Theorem \ref{thm : StS to fields}.

\section{Preliminaries} \label{sec:prelim}

In this section we set up notation, define several bilinear forms that appear in the Euler-Lagrange equations and prove some preliminary results.

\subsection{Algebraic operations}
The inner products for smooth sections give rise to $L^2$-bundle metrics
\begin{align*}
	\langle \alpha, \beta \rangle_{L^2, (S \otimes F)_+} &:= \int_\M \langle \alpha, \beta \rangle_{(\Delta \otimes \V)_+} \,d\vol  && \forall \alpha, \beta \in C^\infty(\M, (\Delta \otimes \V)_+),\\
	\langle \alpha, \beta \rangle_{L^2, \Ad} &:= \int_\M \langle \alpha, \beta \rangle_{\g} \,d\vol  && \forall \alpha, \beta \in C^\infty(\M, \g),\\ 
	\langle \alpha, \beta \rangle_{L^2, E} &:= \int_\M \langle \alpha, \beta \rangle_{\W} \,d\vol  && \forall \alpha, \beta \in C^\infty(\M, \W).
\end{align*}

For what follows, it is convenient to define actions of $X\in \Omega^m(\M, \g) $, $m=0,1,2$, on $\psi\in C^\infty(\M, (\Delta \otimes \V))$.
If we write \begin{align*}
\psi &:= s \otimes \tau \quad \mbox{for some $s \in C^\infty(\M, \Delta)$ and $\tau \in C^\infty(\M, \V)$,}\end{align*}  the actions take the form
\begin{align}\label{eq : explicit bullet}\left\{\begin{aligned}
		X \bullet \psi &:=   s \otimes \left(\varrho_{*}(X) \tau\right) = \varrho_{*}(X) \psi&& X \in C^\infty(\M, \g),\\
		X \bullet \psi &:= \imath \sum_{k=0}^3 \Gamma^k s \otimes \varrho_{*}\left(X(e_k)\right) \tau && X \in \Omega^1(\M, \g),\\
		X \bullet \psi &:= -\sum_{i,j=0}^3 \Gamma^i \Gamma^j s \otimes \varrho_{*}\left(X(e_i, e_j)\right) \tau && X \in \Omega^2(\M, \g).
	\end{aligned}\right.\end{align}
	(Observe that the first product does not require the Clifford structure.)

\subsection{The Euler-Lagrange equations}

Let $(\psi_L, \psi_R, A, \Phi)$ be a critical point of $\mathscr{A}_{\mathrm{SM}}$ in $\mathbb{D}$. In order to derive the Euler-Lagrange equations we consider a variation of the fields around $(\psi_L, \psi_R, A, \Phi)$ given by
\[\begin{aligned}	\psi_{L, t_L} & = \psi_L  + t_L \varphi_L, &&  \forall t_L \geq 0\\
	\psi_{R, t_R} & = \psi_R  + t_R \varphi_R, &&  \forall t_R \geq 0\\
	A_{t_{\mathrm{YM}}} &= A + t_{\mathrm{YM}} W, &&  \forall t_{\mathrm{YM}} \geq 0\\
	\Phi_{t_{\mathrm{H}}} & = \Phi + t_{\mathrm{H}} \Upsilon, &&  \forall t_{\mathrm{H}} \geq 0
.
\end{aligned}\]
Taking the partial derivatives of the resulting variation $$\mathscr{A}_{\mathrm{SM}}[A_{t_{\mathrm{YM}}}, \Phi_{t_{\mathrm{H}}}, \psi_{L, t_L}, \psi_{R, t_R}]$$  in $t_{\mathrm{YM}}$, $t_{\mathrm{H}}$, $t_L$, $t_R$ at $0$ respectively, and using the self-adjointness of $\slashed{D}_A$, yields the Euler-Lagrange equations of the Standard Model,
\begin{align}\label{L}
	\slashed{D}_A \psi_L &=          \J_{\mathrm{YH}, L}(\Phi, \psi_R)
	\\ \label{R}
	\slashed{D}_A \psi_R &=          \J_{\mathrm{YH}, R}(\psi_L, \Phi)
	\\
	\label{YM}
	D_A^\ast F_{A}    &=	\J_{\mathrm{YMH}}(d_A \Phi, \Phi)  + \J_{\mathrm{YMD}}^1(\psi_L, \psi_L) +   \J_{\mathrm{YMD}}^1(\psi_R, \psi_R)  
	\\\label{H}
	d_A^\ast d_A \Phi  &=      \mathbf{V}'(|\Phi|_{\W}^2) \Phi    + \J_{\mathrm{HY}}(\psi_L, \psi_R),
\end{align} where the bilinear interaction forms $\J_\cdot$ are given below.
\begin{itemize}

	\item  The bilinear forms 
	\begin{align*}\J_{\mathrm{YH},L} : C^\infty(\M, \W)\times C^\infty(\M, \Delta_R \otimes \V_{R})&\longrightarrow C^\infty(\M, \Delta_L \otimes \V_{L})\\
		\J_{\mathrm{YH},R} : C^\infty(\M, \Delta_L \otimes \V_{L})\times C^\infty(\M, \W)&\longrightarrow C^\infty(\M, \Delta_R \otimes \V_{R})\end{align*}
	are defined by
	\begin{align*}\langle \J_{\mathrm{YH},L}(\Upsilon,\varphi_R),\varphi_L\rangle_{S\otimes F}&:=-\frac{1}{2}     \mathbf{Y}(\varphi_L,\Upsilon,\varphi_R), \quad \forall \varphi_L \in C^\infty(\M, \Delta_L \otimes \V_L)
	\\
	\langle \J_{\mathrm{YH},R}(\varphi_L,\Upsilon),\varphi_R\rangle_{S\otimes F}&:=-\frac{1}{2}     \mathbf{Y}(\varphi_L,\Upsilon,\varphi_R), \quad \forall \varphi_R \in C^\infty(\M, \Delta_R \otimes \V_R).\end{align*}

	 \item The {\bf real} bilinear form $$\J_{\mathrm{YMH}}:   C^\infty(\M, \W)  \times C^\infty(\M, \W)\longrightarrow C^\infty(\M, \g)$$ is  defined 
	 by	 
	\begin{equation*}
\langle\J_{\mathrm{YMH}}(\Upsilon_1, \Upsilon_2), W\rangle_{\text{\rm Ad}} :=	2	\Re \langle \Upsilon_1,\rho_{*}(W)\Upsilon_2\rangle_{E},\quad \forall W\in C^\infty(\M, \g).
	\end{equation*}

\item Let $k = 0, 1$. The {\bf real} bilinear form $$\J_{\mathrm{YMD}}^k:   C^\infty\left(\M, \Delta \otimes \V\right) \times C^\infty\left(\M, \Delta \otimes \V\right)    \longrightarrow \Omega^k(\M, \g)$$ is  defined 
by 
	\begin{equation*}
\langle\J_{\mathrm{YMD}}^k(\varphi_1, \varphi_2), W\rangle_{\text{\rm Ad}} :=		\Re \langle \varphi_1,   W \bullet \varphi_2\rangle_{(S\otimes F)_+},\quad \forall W\in \Omega^k(\M, \g).
		\label{eq:ex2}
	\end{equation*} 
Only the case $k=1$ appears in the Euler-Lagrange equations, but the case $k=0$ is needed later when formulating the compatibility condition, see Proposition \ref{prop_comp_cond}.

  \item  The bilinear form $\J_{\mathrm{HY}}$   $$\J_{\mathrm{HY}}:   C^\infty(\M, \Delta_L \otimes \V_L)  \times C^\infty(\M, \Delta_R \otimes \V_R)\longrightarrow C^\infty(\M, \W)$$  
 is defined by
 \[  \langle \J_{\mathrm{HY}}(\varphi_L, \varphi_R), \Upsilon \rangle_E :=   -\frac{1}{2}     \mathbf{Y}(\varphi_L, \Upsilon,\varphi_R), \quad \forall \Upsilon \in C^\infty(\M, \W). \]

 \end{itemize}

We make several remarks with regard to the bilinear forms defined above. 
\begin{itemize}
\item The subscript letters in $\J_\cdot$ indicate the Lagrangians that originate the bilinear form.

	\item  {\rm Since the inner product in $\mathcal W$ is invariant under $\rho(g)$ for all $g\in G$ we see that
		$\rho_{*}(X)$ is skew-Hermitian for all $X\in \g$ and thus $\J_{\mathrm{YMH}}$ is anti-symmetric. }

\item {\rm Likewise
		$\varrho_{*}(X)$ is skew-Hermitian for all $X\in \g$. But the Clifford product is also anti-symmetric thus the combination of both produces a form $\J_{\mathrm{YMD}}^1$ that is symmetric unlike $\J_{\mathrm{YMH}}$. }
	
	\item For the sake of convenience, we call  \eqref{L} and \eqref{R} the \textbf{left and right Dirac channels}, \eqref{YM} the \textbf{Yang--Mills channel}, and \eqref{H} the \textbf{Higgs channel}. Their perturbations and linearizations will inherit their name in the perturbed and linearized systems of \eqref{L}-\eqref{H} which will appear below.
 
\end{itemize}

\subsection{The Lichnerowicz-Weitzenb\"ock formula}
In order to convert the first order Dirac channels \eqref{L}-\eqref{R} into a wave equations we will need an expression for $\slashed{D}_A^2$. This is given by a standard Lichnerowicz-Weitzenb\"ock type formula:

 \begin{equation}\label{lichnerowicz}\slashed{D}_A^2 \varphi = 	d^\ast_{  A}d_{A} \varphi + \frac{1}{2} F_A \bullet \varphi,\end{equation}   
(see \cite[(1.1)]{GL_80} for the Riemannian version).  For completeness, we provide a proof of \eqref{lichnerowicz}. We first expand $\slashed{D}_A^2$, by definition,
	\begin{align*}
		\slashed{D}_A^2 \varphi = - \Gamma^b \Gamma^a \left(  \p_b \p_a \varphi  +  (\varrho_\ast A_a) \p_b \varphi + \p_b(\varrho_\ast  A_a) \psi + (\varrho_\ast  A_b) \p_a \varphi + (\varrho_\ast  A_b) (\varrho_\ast  A_a) \varphi  \right).
	\end{align*}
	On one hand, the diagonal terms of $\slashed{D}_A^2 \varphi$ is just   $\Box_A \varphi$, which, by \cite[(16)]{CLOP1}, reads
	\begin{align*}
		d^\ast_{  A}d_{  A} = - \Gamma^b \Gamma^b \left(  \p_b \p_b \varphi + 2 (\varrho_\ast A_b) \p_b \varphi +  \p_b (\varrho_\ast A_b) \varphi +  (\varrho_\ast A_b)  (\varrho_\ast A_b) \varphi  \right).
	\end{align*}
	On the other hand, the off-diagonal terms  can be computed as follows
	\begin{align*}
		\lefteqn{ -   \sum_{b \neq a} \Gamma^b \Gamma^a  \left(  \p_b \p_a \varphi  +  (\varrho_\ast A_a) \p_b \varphi + \p_b (\varrho_\ast A_a) \varphi + (\varrho_\ast A_b) \p_a \varphi + (\varrho_\ast A_b) (\varrho_\ast A_a) \varphi  \right) } 
        \\ 
		&=      \sum_{a < b} \Gamma^b \Gamma^a  \left(  \p_a \p_b \varphi  +  (\varrho_\ast A_b) \p_a \varphi + \p_a (\varrho_\ast A_b) \varphi + (\varrho_\ast A_a) \p_b \varphi + (\varrho_\ast A_a) (\varrho_\ast A_b) \varphi  \right) \\
		&\quad  -   \sum_{a < b} \Gamma^b \Gamma^a  \left(  \p_b \p_a \varphi  +  (\varrho_\ast A_a) \p_b \varphi + \p_b (\varrho_\ast A_a) \varphi + (\varrho_\ast A_b) \p_a \varphi + (\varrho_\ast A_b) (\varrho_\ast A_a) \varphi  \right) \\
		&=       \sum_{a < b} \Gamma^b \Gamma^a \left(     \p_a (\varrho_\ast A_b) \varphi  + (\varrho_\ast A_a) (\varrho_\ast A_b) \varphi   - \p_b (\varrho_\ast A_a) \varphi   - (\varrho_\ast A_b) (\varrho_\ast A_a) \varphi  \right)\\
		&=        \sum_{a < b} \Gamma^b \Gamma^a (\varrho_\ast F_A(e_a, e_b)) \varphi 
		=      \frac{1}{2} F_A \bullet \varphi.
	\end{align*}

\subsection{Dirac action on bilinear forms}

When converting the Dirac equations into wave equations, we will need the following.

\begin{lemma}	\label{lemma:DAJ}\begin{align*}
		\slashed{D}_A \J_{\mathrm{YH}, L}(\Phi, \psi_R) &= \J_{\mathrm{YH}, L}(d_A \Phi(e_a), \imath \Gamma^a \psi_R) + \J_{\mathrm{YH}, L}(\Phi, \J_{\mathrm{YH}, R}(\psi_L, \Phi))\\
		\slashed{D}_A \J_{\mathrm{YH}, R}(\psi_L, \Phi) &= \J_{\mathrm{YH}, R}(\imath \Gamma^a \psi_L, d_A \Phi(e_a)) + \J_{\mathrm{YH}, R}(\J_{\mathrm{YH}, L}(\Phi, \psi_R), \Phi).
\end{align*}\end{lemma} 

\begin{proof} We only give the proof for the first identity since the second is similar.  
	%
	%
	The Yukawa coupling $\mathbf{Y}$, with fixed sections $\Phi \in C^\infty(\M, \W)$, $\psi_{L} \in C^\infty(\M, \Delta_L \otimes \V_L)$ and $\psi_{R} \in C^\infty(\M, \Delta_R \otimes \V_R)$, induces a $1$-form   $\omega$ through
	\[\omega(X):=\mathbf{Y}(X\cdot \psi_{L},\Phi, \psi_{R}), \quad \forall X := X^\alpha e_\alpha\in \mathfrak{X}(\M),\] where the Clifford action $X\cdot \psi_{L}$ is defined by 
	\[ X \cdot  \psi_L := \imath X^\alpha \Gamma_\alpha \psi_L. \]
	If we write $\omega := \omega_\alpha dx^\alpha$ then
	\[\omega^{\alpha}=\mathbf{Y}(\imath \Gamma^{\alpha} \psi_{L},\Phi, \psi_{R}).\] 
	Recall that	\[ d^\ast \omega = \star \,d\,\star\omega=-\partial_{\alpha}\omega^{\alpha}.\]
	Now we compute using that $\mathbf{Y}$ is trilinear
	\begin{align*}
		\partial_{\alpha}\omega^{\alpha}&=\mathbf{Y}(\imath\Gamma^{\alpha}\partial_{\alpha}\psi_{L},\Phi,\psi_{R})+\mathbf{Y}(\imath\Gamma^{\alpha} \psi_{L},\partial_{\alpha}\Phi,\psi_{R})+\mathbf{Y}(\imath\Gamma^{\alpha} \psi_{L},\Phi,\partial_{\alpha}\psi_{R})\\
		&=\mathbf{Y}(\imath\Gamma^{\alpha}\partial_{\alpha}\psi_{L},\Phi,\psi_{R})+\mathbf{Y}(\imath\Gamma^{\alpha} \psi_{L},\partial_{\alpha}\Phi,\psi_{R})-\mathbf{Y}(\psi_{L},\Phi,\imath\Gamma^{\alpha}\partial_{\alpha}\psi_{R}).
	\end{align*}
	For the second equality we used that the Clifford action is antisymmetric with respect to the Dirac bundle metric. 
	We replace $d$ by $\slashed{D}_A$ and $d_{A, \rho}$ and get
	\begin{align*}
		\partial_{\alpha}\omega^{\alpha}&=\mathbf{Y}(\slashed{D}_{A}\psi_{L},\Phi,\psi_{R})-\mathbf{Y}(\imath\Gamma^{\alpha}\varrho_{L*}(A_{\alpha})\psi_{L},\Phi,\psi_{R})\\
		&\quad +\mathbf{Y}(\imath\Gamma^{\alpha} \psi_{L},d_{A, \rho}\Phi(e_{\alpha}),\psi_{R})-\mathbf{Y}(\imath\Gamma^{\alpha} \psi_{L},\rho_{*}(A_{\alpha})\Phi,\psi_{R})\\&\quad -\mathbf{Y}( \psi_{L},\Phi,\slashed{D}_{A}\psi_{R})+\mathbf{Y}(\psi_{L},\Phi,\imath\Gamma^{\alpha}\varrho_{R*}(A_{\alpha})\psi_{R}).
	\end{align*}
	Using antisymmetry again we obtain
	\begin{align}
		\notag	\partial_{\alpha}\omega^{\alpha}&=\mathbf{Y}(\slashed{D}_{A}\psi_{L},\Phi,\psi_{R})  -\mathbf{Y}(\imath\Gamma^{\alpha}\varrho_{L*}(A_{\alpha})\psi_{L},\Phi,\psi_{R})  \\ \notag& \quad +\mathbf{Y}(\imath\Gamma^{\alpha} \psi_{L},d_{A, \rho}\Phi(e_{\alpha}),\psi_{R})-\mathbf{Y}(\imath\Gamma^{\alpha} \psi_{L},\rho_{*}(A_{\alpha})\Phi,\psi_{R})
		\\ \notag 
		&   \quad -\mathbf{Y}(\psi_{L},\Phi,\slashed{D}_{A}\psi_{R})-\mathbf{Y}(\imath\Gamma^{\alpha}\psi_{L},\Phi,\varrho_{R*}(A_{\alpha})\psi_{R}).
	\end{align}
	Combining this with \eqref{eqn : gauge-invariance of Yukawa}    we derive (compare with \cite[Exercise 7.9.12]{Ha}):
	\begin{align*}
		\mathbf{Y}(\slashed{D}_{A}\psi_{L},\Phi,\psi_{R})+\mathbf{Y}(\imath\Gamma^{\alpha} \psi_{L},d_{A}\Phi(e_{\alpha}),\psi_{R})-\mathbf{Y}(\psi_{L},\Phi,\slashed{D}_{A}\psi_{R})=- d^\ast \omega.\end{align*}

	

	
	Using the self-adjointness of $\slashed{D}_A$, the definition of $\J_{\mathrm{YH}, L}(\Phi, \psi_R)$ and integrating $d^\ast \omega$ over $\M$, we obtain
	
	\begin{align*}\lefteqn{\langle  \slashed{D}_A  \J_{\mathrm{YH}, L}(\Phi, \psi_R) ,\psi_{L}\rangle_{S\otimes F, L^{2}}}
		\\&=\langle   \J_{\mathrm{YH}, L}(\Phi, \psi_R) ,\slashed{D}_{A}\psi_{L}\rangle_{S\otimes F,L^2}
		\\&=-\frac{1}{2}\int_{\mathbb{M}}\mathbf{Y}(\slashed{D}_{A}\psi_{L},\Phi,\psi_{R})\,\text{dvol}_{g}
		\\&=-\frac{1}{2}\int_{\mathbb{M}}\mathbf{Y}(\psi_{L},d_{A}\Phi(e_{\alpha}),\imath\Gamma^{\alpha}\psi_{R}),\text{dvol}_{g} - \frac{1}{2}\int_{\mathbb{M}}\mathbf{Y}(\psi_{L},\Phi,\slashed{D}_{A}\psi_{R}),\text{dvol}_{g}.
	\end{align*}
	We conclude the proof by using \eqref{R}.
\end{proof}

\section{Active measurement} \label{sec:sts}

In this section we show how to pass from the data set $\mathcal{D}_{(\psi, A, \Phi)}(\mho)$ to a source-to-solution map $\mathbf{L}_{(\psi, A, \Phi)}$. 
This map can be viewed as a model of measurements, where 
sources acting on $\mho$ produce perturbations to the background fields $(\psi, A, \Phi)$, and the perturbed fields are recorded on $\mho$. 
The correspondence between $\mathbf{L}_{(\psi, A, \Phi)}$ and
$\mathcal{D}_{(\psi, A, \Phi)}(\mho)$ is similar to that between the Dirichlet-to-Neumann map and its graph, that is, the Cauchy data set, in the context of the classical Calder\'on problem.

In the present context, we face several technical complications when defining the source-to-solution map. It is for this reason, that we preferred to formulate our main result using the more abstract set $\mathcal{D}_{(\psi, A, \Phi)}(\mho)$ rather than the more concrete map $\mathbf{L}_{(\psi, A, \Phi)}$. All these complications were already present in our previous works \cite{CLOP2, CLOP3}. Here we omit proofs that only require routine modifications, and focus on the details that differ from \cite{CLOP2, CLOP3} more substantially. 

With a suitable gauge choice, that is, the relative Lorenz gauge, the Euler-Lagrange equations can be reduced to a non-linear hyperbolic system. Due to the compactness of the causal diamond $\mathbb D$, this system has a unique solution when the background solution $(\psi, A, \Phi)$ is perturbed by a small and smooth enough source term (Proposition \ref{prop : solution in the relative gauge}). However, such a solution satisfies the original Euler-Lagrange equations only if it is in the correct gauge. It turns out (Proposition \ref{prop_gauge_comp}), that this is guaranteed by the source satisfying a compatibility condition, that arises naturally from the gauge invariance (Proposition \ref{prop_comp_cond}).

While the relative Lorenz gauge yields a good theory of the direct problem for the Standard Model, we need a second gauge, the temporal gauge, to extract the source-to-solution map from the data set $\mathcal{D}_{(\psi, A, \Phi)}(\mho)$. The natural way to pass to the relative Lorenz gauge is to solve the direct problem as sketched above, but this requires considering the whole diamond $\mathbb D$. On the other hand, the temporal gauge localizes very nicely in $\mho$ (Proposition \ref{prop : uniqueness temporal}).

In summary, the evaluation $\mathbf{L}_{(\psi, A, \Phi)} f$ for a small and smooth enough source $f$, satisfying the compatibility condition, is constructed by first solving the Euler-Lagrange equations in the relative Lorenz gauge and then transforming the solution into the temporal gauge. Due to the localized nature of the temporal gauge, we can relate the restriction of the transformed solution to a point in the data set $\mathcal{D}_{(\psi, A, \Phi)}(\mho)$.

\subsection{Compatibility condition} \label{subsection : compatibility}
We perturb \eqref{L}-\eqref{H} with sources
\[\left(\mathcal{K}, \mathcal{J}, \mathcal{F}\right) \in C_0^\infty(\mho, (\Delta \otimes \V)_+)  \times C_0^\infty(\mho, T^\ast \M \otimes \g) \times C_0^\infty(\mho, \W)  .\]
Consider fields \[\left(\phi, V, \Psi\right) \in  C^\infty(\mathbb{D}, (\Delta \otimes \V)_+) \times C^\infty(\mathbb{D}, T^\ast \M \otimes \g) \times C^\infty(\mathbb{D}, \W)\] solving the perturbed system
\begin{align}	
	\label{pL}
	\slashed{D}_V \phi_L - \J_{\mathrm{YH}, L}(\Psi, \phi_R)&=     \mathcal{K}_L      
	\\
	\label{pR}
	\slashed{D}_V \phi_R - \J_{\mathrm{YH}, R}(\phi_L, \Psi) &=          \mathcal{K}_R
	\\
	\label{pYM}
	D_V^\ast F_{V}    -	\J_{\mathrm{YMH}}(d_V \Psi, \Psi)  - \J_{\mathrm{YMD}}^1(\phi_L, \phi_L) -   \J_{\mathrm{YMD}}^1(\phi_R, \phi_R)  &= \mathcal{J}
	\\
	\label{pH}
	d_V^\ast d_V \Psi   -     \mathbf{V}'(|\Psi|^2_{\W}) \Psi  - \J_{\mathrm{HY}}(\phi_L, \phi_R)  &= \mathcal{F}\end{align} with $\phi = \phi_L + \phi_R$ and $\mathcal{K} = \mathcal{K}_L + \mathcal{K}_R$.

The sources cannot be arbitrarily chosen, and the constraint on the sources is known as the {\it compatibility condition}.
\begin{proposition}\label{prop_comp_cond}
The sources $(\mathcal{K}, \mathcal{J}, \mathcal{F})$ in \eqref{pL}-\eqref{pH} must satisfy the compatibility condition   
	\begin{align}
		\label{eqn : compatibility} D_V^\ast \mathcal{J} +   \J_{\mathrm{YMH}}(\mathcal{F}, \Psi) + 2\J_{\mathrm{YMD}}^0(\mathcal{K}_L, \phi_L) + 2\J_{\mathrm{YMD}}^0(\mathcal{K}_R, \phi_R) = 0.
	\end{align}

\end{proposition}

\begin{proof} Equation \eqref{eqn : compatibility} is a consequence of  the gauge invariance of the action functional $\mathscr{A}_{\mathrm{SM}}[(\phi, V, \Psi)]$. To see this, we let $\mathbf{U}_t$ be a smooth family of gauge transformations, with parameter $t$, and denote $\dot{\mathbf{U}}= \p_t \mathbf{U}|_{t=0}$. Then letting $\mathbf{U}_t$ act on the triple $(\phi, V, \Psi)$ and differentiating at $t=0$ gives 
\begin{align*}  
	\p_t \left.\left(\phi \cdot \mathbf{U}_t \right)\right|_{t=0} &= - \lambda \otimes \varrho_\ast (\dot{\mathbf{U}}) v\\
	\p_t \left.\left(V \cdot \mathbf{U}_t\right)\right|_{t=0} &= d\dot{\mathbf{U}} + [V, \dot{\mathbf{U}}] = D_V \dot{\mathbf{U}} \\ 
	\p_t \left.\left(\Psi \cdot \mathbf{U}_t \right)\right|_{t=0} &= - \rho_\ast (\dot{\mathbf{U}}) \Psi.
\end{align*}
The action functional $\mathscr{A}_{\mathrm{SM}}[(\phi, V, \Psi)]$ is gauge invariant, i.e. 
\[\mathscr{A}_{\mathrm{SM}}[(\phi, V, \Psi)] = \mathscr{A}_{\mathrm{SM}}[(\phi, V, \Psi) \cdot \mathbf{U}_t]\] and thus
\[\p_t \left. \mathscr{A}_{\mathrm{SM}}[(\phi, V, \Psi) \cdot \mathbf{U}_t] \right|_{t=0} \equiv 0.\] 
The chain rule yields that 
\begin{align*}
	\lefteqn{\p_t \left. \mathscr{A}_{\mathrm{SM}}[(\phi, V, \Psi) \cdot \mathbf{U}_t] \right|_{t=0}}  \\   
	&=
		2\Re\left(	\slashed{D}_V \phi_L - \J_{\mathrm{YH}, L}(\Psi, \phi_R), - s_L \otimes \varrho_\ast (\dot{\mathbf{U}}) \tau_L\right)_{L^2, S \otimes F}\\
	&\quad + 2\Re\left(	\slashed{D}_V \phi_R - \J_{\mathrm{YH}, R}(\phi_L, \Psi), - s_R \otimes \varrho_\ast (\dot{\mathbf{U}}) \tau_R\right)_{L^2, S \otimes F} 
	\\
	&\quad- \left( 	D_V^\ast F_{V}    -	\J_{\mathrm{YMH}}(d_V \Psi, \Psi)  - \J_{\mathrm{YMD}, L}(\phi_L, \phi_L) -   \J_{\mathrm{YMD}, R}(\phi_R, \phi_R), D_V \dot{\mathbf{U}}    \right)_{L^2, \Ad} \\
	&\quad + 2\Re\left(  d_V^\ast d_V \Psi   -    \mathbf{V}'(|\Psi|^2_{\W}) \Psi  - \J_{\mathrm{HY}}(\phi_L, \phi_R), - \rho_\ast(\dot{\mathbf{U}}) \Psi \right)_{L^2, E} .	
\end{align*} Inserting the perturbed equations \eqref{pL}-\eqref{pH} and using the definitions of the bilinear forms gives that 
\begin{align*}	\lefteqn{\p_t \left. \mathscr{A}_{\mathrm{SM}}[(\phi, V, \Psi) \cdot \mathbf{U}_t] \right|_{t=0}}  \\
	&=	 2\Re\left(	\mathcal{K}_L, - s_L \otimes \varrho_\ast (\dot{\mathbf{U}}) \tau_L\right)_{L^2, S \otimes F} 
  + 2\Re\left(	\mathcal{K}_R, - s_R \otimes \varrho_\ast (\dot{\mathbf{U}}) \tau_R\right)_{L^2, S \otimes F} \\ &\quad -\left( 	\mathcal{J}, D_V \dot{\mathbf{U}}    \right)_{L^2, \Ad}  
   + 2\Re\left(  \mathcal{F}, - \rho_\ast(\dot{\mathbf{U}}) \Psi \right)_{L^2, E} \\
	&= 	 - 2\left(	\J_{\mathrm{YMD}}^0(\mathcal{K}_L, \psi_L),  \dot{\mathbf{U}} \right)_{L^2, S \otimes F}
 - 2\left(	\J_{\mathrm{YMD}}^0(\mathcal{K}_R, \psi_R),  \dot{\mathbf{U}} \right)_{L^2, S \otimes F}
	\\&\quad-\left(  D_V^\ast	\mathcal{J}, \dot{\mathbf{U}}    \right)_{L^2, \Ad} 
	 -  \left( \J_{\mathrm{YMH}}(\mathcal{F}, \Psi),   \dot{\mathbf{U}}  \right)_{L^2, E} 
.	
\end{align*} 
This proves \eqref{eqn : compatibility}.
\end{proof}

\subsection{Equations for the perturbations} \label{subsec : hyperbolic}

We will next rewrite the system \eqref{pL}-\eqref{pH} in terms of the perturbations 
\begin{align*}
	\varphi_{L} := \phi_{L} - \psi_{L}, \quad \varphi_{R} := \phi_{R} - \psi_{R},
\quad
W := V - A, \quad
	\Upsilon := \Psi - \Phi.
\end{align*}
This is convenient since the resulting system for $(\varphi, W, \Upsilon)$
can be studied via linearizing it at zero. 
By using  \[F_V = F_A + dW + [A, W] + [W, W]/2,\] we obtain 
the following equations for the perturbations	
\begin{align} 
	\label{eq:perturbed D+ prelim}
	\slashed{D}_A \varphi_L +   W \bullet (\varphi_L + \psi_L) - \mathcal{R}_{\Phi, \psi_R, L}(\Upsilon, \varphi_R)&=     \mathcal{K}_L      
	\\
	\label{eq:perturbed D- prelim}
	\slashed{D}_A \varphi_R +   W \bullet (\varphi_R + \psi_R)  - \mathcal{R}_{\Phi, \psi_L, R}(\Upsilon, \varphi_L) &=     \mathcal{K}_R 
	\\
	D_{A}^*D_{A}W  + \mathcal{N}_{A}(W)+\sum_{i=1}^{3}\mathcal{M}_{A,\Phi}^{i}(W,\Upsilon) -  \mathcal{P}_{\psi_L, \psi_R}(\varphi_L, \varphi_R)     &= \mathcal{J}
	\label{eq:perturbed YM prelim} \\	
	\Box_{A} \Upsilon+ \sum_{i=1}^{3} 	\mathcal{O}_{A, \Phi}^i(W, \Upsilon)  -  \mathcal{Q}_{\psi_L, \psi_R}(\varphi_L, \varphi_R)   &=\mathcal{F},
	\label{eq:perturbed H prelim} 
\end{align} 
where $\Box_{A} = d^\ast_{A} d_{A}$, and we use the shorthand notations
\begin{align*}
	\mathcal{R}_{\Phi, \psi_R, L}(\Upsilon, \varphi_R) &:= \J_{\mathrm{YH}, L}(\Upsilon, \varphi_R) + \J_{\mathrm{YH}, L}(\Phi, \varphi_R) + \J_{\mathrm{YH}, L}(\Upsilon, \psi_R) \\
	\mathcal{R}_{\Phi, \psi_L, R}(\Upsilon, \varphi_L) &:= \J_{\mathrm{YH}, R}(\Upsilon, \varphi_L) + \J_{\mathrm{YH}, R}(\Phi, \varphi_L) + \J_{\mathrm{YH}, R}(\Upsilon, \psi_L)
	\end{align*}
in the Dirac channels,	and
	\begin{align*}
	 \mathcal N_{A}(W) &:=\star[W,\star F_{A}]+\frac{1}{2} D_A^\ast [W, W] + \star [W, \star D_A W] + \frac{1}{2} \star [W, \star [W, W]] \\
	\mathcal{M}_{A,\Phi}^{1}(W,\Upsilon)	&:=\J_{\mathrm{YMH}}(d_{A}\Upsilon,\Phi)+\J_{\mathrm{YMH}}(d_{A}\Phi,\Upsilon)+\J_{\mathrm{YMH}}(\rho_{*}(W)\Phi,\Phi)\\
	\mathcal{M}_{A,\Phi}^{2}(W,\Upsilon) 	&:=\J_{\mathrm{YMH}}(d_{A}\Upsilon,\Upsilon)+\J_{\mathrm{YMH}}(\rho_{*}(W)\Upsilon,\Phi)+\J_{\mathrm{YMH}}(\rho_{*}(W)\Phi,\Upsilon)\\
	\mathcal{M}_{A,\Phi}^{3}(W,\Upsilon) 	&:=\J_{\mathrm{YMH}}(\rho_{*}(W)\Upsilon,\Upsilon)\\
	\mathcal{P}_{\psi_L, \psi_R}(\varphi_L, \varphi_R)  &:= \J_{\mathrm{YMD}, L}(\psi_L, \varphi_L) + \J_{\mathrm{YMD}, L}(\varphi_L, \psi_L) + \J_{\mathrm{YMD}, L}(\varphi_L, \varphi_L) \\ &\quad +  \J_{\mathrm{YMD}, R}(\psi_R, \varphi_R) +  \J_{\mathrm{YMD}, R}(\varphi_R, \psi_R)
   +  \J_{\mathrm{YMD}, R}(\varphi_R, \varphi_R)\\ 
	\mathcal{O}_{A, \Phi}^1(W, \Upsilon) &:= 
	d_{A}^*(\rho_{*}(W)\Phi)+\star(\rho_{*}(W)\wedge\star d_{A}\Phi)
	+ 2 \Re \langle \Phi, \Upsilon \rangle_{\mathcal W} \Phi +  |\Phi|^2_{\W} \Upsilon
	\\
	\mathcal{O}_{A, \Phi}^2(W, \Upsilon) &:= \star(\rho_{*}(W)\wedge\star \rho_{*}(W)\Phi) 
	+ 
	d_{A}^*(\rho_{*}(W)\Upsilon)+\star(\rho_{*}(W)\wedge\star d_{A}\Upsilon)
	\\&\qquad + 2 \Re \langle \Phi, \Upsilon \rangle_{\mathcal W} \Upsilon + |\Upsilon|^2_{\W} \Phi
	\\
	\mathcal{O}_{A, \Phi}^3(W, \Upsilon) &:= \star(\rho_{*}(W)\wedge\star \rho_{*}(W)\Upsilon) 
	+ |\Upsilon|^2_{\W} \Upsilon 	\\
	\mathcal{Q}_{\psi_L, \psi_R}(\varphi_L, \varphi_R)  &:= \J_{\mathrm{HY}}(\psi_L, \varphi_R) + \J_{\mathrm{HY}}(\varphi_L, \psi_R) + \J_{\mathrm{HY}}(\varphi_L, \varphi_R)
\end{align*} in the Yang--Mills and Higgs channels.
We have deviated from the notation used in \cite{CLOP3}, where the linear term $\star[W,\star F_{A}]$ is not a part of $\mathcal N_A(W)$. This term is negligible from the point of view of the method we use.

\subsection{Passing from Dirac to wave equations}
\label{sec_from_Dirac_to_wave}

Applying the Dirac operator to the Dirac channels \eqref{eq:perturbed D+ prelim}-\eqref{eq:perturbed D- prelim} turns them to wave equations. 
In this section we show that the resulting equations are equivalent to the original Dirac equations.  
The rationale behind this transformation is that it allows us to treat all the four equations in a unified manner: the Higgs channel \eqref{eq:perturbed H prelim} is a wave equation and so is the Yang-Mills channel \eqref{eq:perturbed YM prelim} in the relative Lorentz gauge, see \eqref{eq:perturbed YM rL} below.

We apply $\slashed{D}_{A + W}$ to \eqref{eq:perturbed D+ prelim}-\eqref{eq:perturbed D- prelim} and use the Lichnerowicz-Weitzenb\"ock formula \eqref{lichnerowicz} to get
\begin{align} 
	\label{eq:perturbed D+ convt}
	\Box_A \varphi_L + \frac{1}{2} D_A W \bullet \psi_L  + \sum_{i=1}^3\mathcal{T}^i_{A, \psi_L}(W, \varphi_L)   - \slashed{D}_{A + W}\mathcal{R}_{\Phi, \psi_R, L}(\Upsilon, \varphi_R)&=     \tilde{\mathcal{K}}_L      
	\\
	\label{eq:perturbed D- convt}
	\Box_A \varphi_R + \frac{1}{2} D_A W \bullet \psi_R  + \sum_{i=1}^3\mathcal{T}^i_{A, \psi_R}(W, \varphi_R)  - \slashed{D}_{A + W}\mathcal{R}_{\Phi, \psi_L, R}(\Upsilon, \varphi_L) &=    \tilde{\mathcal{K}}_R 
\end{align}
where we introduce additional shorthand notations
\begin{align*}  
	\mathcal{T}^1_{A, \psi}(W, \varphi) &:=  \frac{1}{2} F_A \bullet \varphi + 2 \star (\varrho_\ast(W) \wedge \star d_{A} \psi) \\
	\mathcal{T}^2_{A, \psi}(W, \varphi) &:=\frac{1}{2} D_A W \bullet  \varphi  + 2 \star (\varrho_\ast(W) \wedge \star d_{A} \varphi) \\ 
	& \quad \quad + \star (\varrho_\ast(W) \wedge \star (\varrho_\ast(W) \psi)) + \frac{1}{4} [W, W] \bullet \psi\\
	\mathcal{T}^3_{A, \psi}(W, \varphi) &:=  \star (\varrho_\ast(W) \wedge \star (\varrho_\ast(W) \varphi)) + \frac{1}{4} [W, W] \bullet \varphi,
\end{align*}
and 
	\begin{align}\label{def_K_tilde}
	\tilde{\mathcal{K}}_{L} := \slashed{D}_{A + W}\mathcal{K}_{L},
\quad
	\tilde{\mathcal{K}}_{R} := \slashed{D}_{A + W}\mathcal{K}_{R}.
	\end{align}

For $\tilde{\mathcal{K}}_{L}$ and $\tilde{\mathcal{K}}_{R}$
of the form \eqref{def_K_tilde}, \eqref{eq:perturbed D+ convt}-\eqref{eq:perturbed D- convt} are indeed equivalent to \eqref{eq:perturbed D+ prelim}-\eqref{eq:perturbed D- prelim}.  It suffices to prove that the solution to \eqref{eq:perturbed D+ convt}-\eqref{eq:perturbed D- convt} also solves \eqref{eq:perturbed D+ prelim}-\eqref{eq:perturbed D- prelim}. Writing 
\begin{align*}
	r_L &:=	\slashed{D}_A \varphi_L +   W \bullet (\varphi_L + \psi_L) - \mathcal{R}_{\Phi, \psi_R, L}(\Upsilon, \varphi_R) -    \mathcal{K}_L      
	\\
	r_R &:= \slashed{D}_A \varphi_R +   W \bullet (\varphi_R + \psi_R)  - \mathcal{R}_{\Phi, \psi_L, R}(\Upsilon, \varphi_L) -     \mathcal{K}_R,
\end{align*} 
equations \eqref{eq:perturbed D+ prelim}-\eqref{eq:perturbed D- prelim} are equivalent to $r_L$ and $r_R$ vanishing. 
Applying $\slashed{D}_{A+W}^2$ to $r_L$ and $r_R$ respectively and then using \eqref{lichnerowicz} leads to a symmetric hyperbolic system
\begin{align*}
	d^\ast_{A + W} d_{A + W} r_L + \frac{1}{2} F_{A + W} \bullet r_L &= 0\\
	d^\ast_{A + W} d_{A + W} r_R + \frac{1}{2} F_{A + W} \bullet r_R &= 0.
\end{align*}
If the perturbations $(\varphi, W, \Upsilon)$ vanish near $\p^- \mathbb D$, then $r_L$ and $r_R$ also vanish there, and it follows from the routine local uniqueness theorem for symmetric hyperbolic systems,
that $r_L$ and $r_R$ vanish on the whole $\mathbb D$.

\subsection{Relative Lorenz gauge}\label{subsec : Lorenz}

The relative Lorenz gauge condition
\[D_A^\ast W = 0\]
turns equations \eqref{eq:perturbed YM prelim}-\eqref{eq:perturbed D- convt} to the system of wave equations
\begin{align}
	\label{eq:perturbed D+ rL}
	\Box_{A, \varrho} \varphi_L + \frac{1}{2} D_A W \bullet \psi_L  + \sum_{i=1}^3\mathcal{T}^i_{A, \psi_L}(W, \varphi_L)   - \slashed{D}_{A + W}\mathcal{R}_{\Phi, \psi_R, L}(\Upsilon, \varphi_R)&=     \tilde{\mathcal{K}}_L      
	\\
	\label{eq:perturbed D- rL}
	\Box_{A, \varrho} \varphi_R + \frac{1}{2} D_A W \bullet \psi_R  + \sum_{i=1}^3\mathcal{T}^i_{A, \psi_R}(W, \varphi_R)  - \slashed{D}_{A + W}\mathcal{R}_{\Phi, \psi_L, R}(\Upsilon, \varphi_L) &=     \tilde{\mathcal{K}}_R 
	\\
	\Box_{A, \Ad} W +\mathcal{N}_{A}(W)+\sum_{i=1}^{3}\mathcal{M}_{A,\Phi}^{i}(W,\Upsilon) -  \mathcal{P}_{\psi_L, \psi_R}(\varphi_L, \varphi_R)     &=\mathcal{J}
	\label{eq:perturbed YM rL} \\	
\Box_{A, \rho} \Upsilon+ \sum_{i=1}^{3} 	\mathcal{O}_{A, \Phi}^i(W, \Upsilon)  -  \mathcal{Q}_{\psi_L, \psi_R}(\varphi_L, \varphi_R)   &=\mathcal{F},
	\label{eq:perturbed H rL} 
\end{align}
since  
\begin{align*}
	\Box_{A, \Ad} W &= D_A^\ast D_A  W + D_A D_A^\ast W && \forall W \in C^\infty(\M, T^\ast \M \otimes \g).
\end{align*}
Here the notation emphasizes the fact that, although both the wave operators $\Box_{A, \varrho}$ and $\Box_{A, \rho}$ have the same abstract form $d^\ast_A d_A$, the respective covariant derivatives $d_A$, as defined in Section \ref{sec_diff_ops}, are associated to the distinct representations $\varrho$ and $\rho$.

We proceed to rewrite the compatibility condition \eqref{eqn : compatibility} as an ordinary differential equation for $\mathcal J_0$, 
 \begin{multline}\label{eqn : compatibility bg pt} 
		\p_0 \mathcal{J}_0 + [A_0, \mathcal{J}_0] + [W_0, \mathcal{J}_0] = \p^j \mathcal{J}_j + [A^j, \mathcal{J}_j]  + [W^j, \mathcal{J}_j] \\ 
		- \J_{\mathrm{YMH}}(\mathcal{F}, \Phi + \Upsilon) - 2\J_{\mathrm{YMD}}^0(\mathcal{K}_L, \psi_L + \varphi_L) - 2\J_{\mathrm{YMD}}^0(\mathcal{K}_R, \psi_R + \varphi_R).
	\end{multline} 
Recall that $\tilde{\mathcal{K}}_{L}$ and $\tilde{\mathcal{K}}_{R}$
should be of the form \eqref{def_K_tilde}.
The next two propositions say that \eqref{def_K_tilde}-\eqref{eqn : compatibility bg pt} form a well-posed system, equivalent to the original Euler-Lagrange equations \eqref{pL}-\eqref{pH} with sources.

\begin{proposition}\label{prop : solution in the relative gauge} 
Let $(\psi, A, \Phi)$ be a smooth critical point of $\mathscr{A}_{\mathrm{SM}}[\psi, A, \Phi]$ and let $k \ge 4$. Then for any sufficiently small sources 
	\begin{align}\label{def_source_direct}
\left(\mathcal{K}, \mathcal{J}_1, \mathcal{J}_2, \mathcal{J}_3, \mathcal{F}\right)  \in  H_0^{k+2} \left(  \mho,   \Delta \otimes \V  \oplus \g \oplus \g \oplus \g \oplus \W \right)
	\end{align}
the system \eqref{def_K_tilde}-\eqref{eqn : compatibility bg pt} has a unique solution 
\[ \left(\varphi, W, \Upsilon, \mathcal{J}_0\right)  \in  H^{k+1} \left(  \mathbb D, T^\ast \mathbb D \otimes \g \oplus \W \oplus \Delta \otimes \V  \oplus \g \right).\] 
vanishing near $\p^- \mathbb D$.
The solution depends on $\left(\mathcal{J}_1, \mathcal{J}_2, \mathcal{J}_3, \mathcal{F}, \mathcal{K}\right)$ smoothly. 
\end{proposition} 

To prove this proposition, one can adapt the arguments for \cite[Proposition 3]{CLOP3}.

\begin{proposition}\label{prop_gauge_comp}
Let $(\psi, A, \Phi)$ be a smooth critical point of $\mathscr{A}_{\mathrm{SM}}[\psi, A, \Phi]$, and let $(\varphi, W, \Upsilon, \mathcal J_0)$ vanish near $\p^- \mathbb D$ and solve the system \eqref{def_K_tilde}-\eqref{eqn : compatibility bg pt}
with a source \eqref{def_source_direct}. Then
the perturbed fields $$(\phi, V, \Psi) = (\psi, A, \Phi) + (\varphi, W, \Upsilon)$$ solve \eqref{pL}-\eqref{pH}.
\end{proposition} 
\begin{proof}
Observe that \eqref{eq:perturbed H rL} is \eqref{eq:perturbed H prelim}, and \eqref{eq:perturbed D+ rL}-\eqref{eq:perturbed D- rL} are \eqref{eq:perturbed D+ convt}-\eqref{eq:perturbed D- convt}.
Equations \eqref{eq:perturbed D+ prelim}-\eqref{eq:perturbed D- prelim} follow from \eqref{eq:perturbed D+ convt}-\eqref{eq:perturbed D- convt} and \eqref{def_K_tilde}, together with the vanishing of $(\varphi, W, \Upsilon, \mathcal J_0)$ near $\p^- \mathbb D$, as explained in Section \ref{sec_from_Dirac_to_wave}.
Recalling that \eqref{eq:perturbed D+ prelim}-\eqref{eq:perturbed H prelim} are 
\eqref{pL}-\eqref{pH}, written in terms of perturbations rather than total fields, it remains to show \eqref{pYM}.
Furthermore, as \eqref{eq:perturbed YM rL} and \eqref{eq:perturbed YM prelim} differ only by the term $D_A D_A^\ast W$, and as \eqref{eq:perturbed YM prelim} is \eqref{pYM}, written in terms of perturbations, we have
	\begin{align}\label{gauge_comp_aux}
D_A D_A^\ast W +	D_V^\ast F_{V}    
-\J_{\mathrm{YMH}}(d_V \Psi, \Psi)  
-\J_{\mathrm{YMD}}^1(\phi_L, \phi_L) 
-\J_{\mathrm{YMD}}^1(\phi_R, \phi_R)   
= \mathcal{J}.
	\end{align}
Hence it is enough to show that $D_A^\ast W = 0$.
By applying $D_V^\ast$ to the above equation we will show that
	\begin{align}\label{wave_eq_div_W}
D_V^* D_A D_A^\ast W = 0.
	\end{align} 
This is a linear wave equation for $D_A^\ast W$. As $W$ vanishes near $\p^- \mathbb D$, we can conclude that $D_A^\ast W$ vanishes in $\mathbb D$.

It remains to show \eqref{wave_eq_div_W}.
We begin by considering the coupling terms in \eqref{gauge_comp_aux}. For any $X \in C_0^\infty(\M,   \g)$, we check that 
\begin{align*}
&\langle   D_V^\ast \J_{\mathrm{YMH}}(d_V \Psi, \Psi) , X \rangle_{\Ad, L^2} 
= 2 \Re \langle d_V \Psi, \rho_\ast(D_V X) \Psi \rangle_{E, L^2}
\\&\qquad= 2 \Re \langle d_V \Psi, (d_V\rho_\ast( X)) \Psi \rangle_{E, L^2}
= 2 \Re \langle d_V \Psi, d_V(\rho_\ast( X) \Psi) -  \rho_\ast( X) d_V\Psi  \rangle_{E, L^2}
\\&\qquad= 2\Re \langle  d^\ast_V d_V \Psi, \rho_\ast (X) \Psi   \rangle_{E, L^2} - 2\Re \langle d_V \Psi, \rho_\ast(X) (d_V \Psi) \rangle_{E, L^2} 
\\&\qquad= 2\Re \langle  d^\ast_V d_V \Psi, \rho_\ast (X) \Psi   \rangle_{E, L^2}.
\end{align*} 
In the last identity, we used the fact that $\rho_\ast (X)$ is skew-Hermitian. Furthermore, 
\begin{align*} 
&\langle   D_V^\ast \J_{\mathrm{YMD}}^1(\phi, \phi), X \rangle_{\Ad, L^2} 
= \Re \langle \phi, (D_V X) \bullet \phi \rangle_{S \otimes F, L^2}
\\&\qquad= \Re \langle \phi,    (\slashed{D}_V (X \bullet \phi) - X \bullet \slashed{D}_V\phi) \rangle_{S \otimes F, L^2}
\\&\qquad= \Re \langle\slashed{D}_V \phi, X \bullet \phi \rangle_{S \otimes F, L^2} - \Re \langle \phi, X \bullet \slashed{D}_V \phi \rangle_{S \otimes F, L^2},
\\&\qquad= 2 \Re \langle\slashed{D}_V \phi, X \bullet \phi \rangle_{S \otimes F, L^2}.
\end{align*} 
In the last identity, we used the fact that $\varrho_*(X)$ is skew-Hermitian. Recall that, as $X \in C_0^\infty(M, \g)$, 
the product $X \bullet \psi$ reads simply $\varrho_*(X) \psi$ for sections $\psi$ of the twisted spinor bundle.

Together with \eqref{pL}-\eqref{pR} and \eqref{pH} the previous two identities imply
\begin{align*}
		\langle   D_V^\ast \J_{\mathrm{YMH}}(d_V \Psi, \Psi) , X \rangle_{\Ad, L^2} &= 2\Re \langle  \mathcal{F}, \rho_\ast (X) \Psi   \rangle_{E, L^2} + 2\Re \langle  \J_{\mathrm{HY}}(\phi_L, \phi_R), \rho_\ast (X) \Psi   \rangle_{E, L^2}   \\
	&=    \langle  \J_{\mathrm{YMH}}(\mathcal{F}, \Psi), X   \rangle_{\Ad, L^2} 
	- \mathbf{Y}(\varphi_L, \rho_\ast(X)\Psi, \varphi_R) 
\\
	\langle   D_V^\ast \J_{\mathrm{YMD}}^1(\phi_L, \phi_L), X \rangle_{\Ad, L^2} &=   2\Re \langle     \mathcal{K}_L , X \bullet \phi_L \rangle_{S \otimes F, L^2} + 2\Re \langle\J_{\mathrm{YH}, L}(\Psi, \phi_R)   , X \bullet \phi_L \rangle_{S \otimes F, L^2}   \\
                                                                                 &=   2\langle\J_{\mathrm{YMD}}^0  (      \mathcal{K}_L, \phi_L) , X  \rangle_{\Ad, L^2}    -  \mathbf{Y}(X \bullet \varphi_L, \Psi, \varphi_R) \\
	\langle   D_V^\ast \J_{\mathrm{YMD}}^1(\phi_R, \phi_R), X \rangle_{\Ad, L^2} 
	&=   2\langle\J_{\mathrm{YMD}}^0  (     \mathcal{K}_R, \phi_R) , X  \rangle_{\Ad, L^2}  -  \mathbf{Y}(\varphi_L,  \Psi, X \bullet \varphi_R),
\end{align*}where we used again the fact that $\rho_*(X)$ is skew-Hermitian to eliminate the Higgs potential term in \eqref{pH}.

Now we assert that \begin{align}\label{eqn : gauge-invariance of Yukawa}\mathbf{Y}(\varphi_L, \rho_\ast(X)\Psi, \varphi_R)  +  \mathbf{Y}(X \bullet \varphi_L, \Psi, \varphi_R)  +  \mathbf{Y}(\varphi_L,  \Psi, X \bullet \varphi_R) =  0.\end{align} This is a direct consequence of the gauge invariance of the Yukawa Lagrangian $\mathbf{Y}$.  Consider a smooth family $\mathbf{U}_t$ of gauges and write  $X = \p_t \mathbf{U} |_{t=0}$. The gauge invariance implies \[\p_t \mathbf{Y}(\phi_L \cdot \mathbf{U}_t, \Psi \cdot \mathbf{U}_t, \phi_R \cdot \mathbf{U}_t)  = 0,\]
and \eqref{eqn : gauge-invariance of Yukawa} follows from the multilinearity of the Yukawa form.

Combining \eqref{eqn : gauge-invariance of Yukawa} with the three equations preceding it gives
	\begin{align*}
&D_V^\ast \J_{\mathrm{YMH}}(d_V \Psi, \Psi)  
+ D_V^\ast \J_{\mathrm{YMD}}^1(\phi_L, \phi_L) 
+ D_V^\ast \J_{\mathrm{YMD}}^1(\phi_R, \phi_R)   
\\&\qquad=
\J_{\mathrm{YMH}}(\mathcal{F}, \Psi)
+ 2 \J_{\mathrm{YMD}}^0(\mathcal{K}_L, \phi_L)
+ 2 \J_{\mathrm{YMD}}^0(\mathcal{K}_R, \phi_R).
	\end{align*}
Recalling that \eqref{eqn : compatibility bg pt} is a rearrangement of the compatibility condition \eqref{eqn : compatibility}, we get 
\begin{align*}
&- D_V^\ast \mathcal{J}  
= 
\J_{\mathrm{YMH}}(\mathcal{F}, \Psi)
+ 2 \J_{\mathrm{YMD}}^0(\mathcal{K}_L, \phi_L)
+ 2 \J_{\mathrm{YMD}}^0(\mathcal{K}_R, \phi_R)
\\&\qquad=
D_V^\ast \J_{\mathrm{YMH}}(d_V \Psi, \Psi) 
+ D_V^\ast \J_{\mathrm{YMD}}^1(\phi_L, \phi_L) 
+ D_V^\ast \J_{\mathrm{YMD}}^1(\phi_R, \phi_R).  
\end{align*}
Equation \eqref{wave_eq_div_W} follows now by applying $D_V^*$ to \eqref{gauge_comp_aux} and using $D_V^\ast D_V^\ast F_{V} = 0$, see e.g. \cite[Lemma 2]{CLOP2}.
\end{proof}

\subsection{Temporal gauge} \label{subsec : temporal}

A connection $V = V_\alpha dx^\alpha$ is said to be in the temporal gauge if $V_0 = 0$. 
The gauge transformation 
	$\mathscr{T} (V) := V \cdot \mathbf{U}$ with $\mathbf{U} \in G^0(\mathbb{D}, p)$ satisfying
\begin{align*}
	\p_{0} \mathbf{U} &= - V_0 \mathbf{U} \\  \mathbf{U}|_{\p^-\mathbb{D}} &= \id
\end{align*}
takes a connection $V \in \Omega^1(\mathbb{D}; \g)$ to the temporal gauge. We write also 
	\begin{align*}
\mathscr{T}(\varphi, W,  \Upsilon) = (\varphi, W,  \Upsilon) \cdot \mathbf{U}.
	\end{align*}

We have the following uniqueness result for \eqref{eq:perturbed YM prelim}-\eqref{eq:perturbed D- convt} in the temporal gauge, with essentially the same proof as that of \cite[Proposition 2]{CLOP2}.

\begin{proposition}\label{prop : uniqueness temporal}   Given a source $(\mathcal{K}, \mathcal{J}, \mathcal{F})$, suppose that
\[\left(V^{(j)}, \Psi^{(j)}, \phi^{(j)}\right) \in C^3 \left(  \mathbb{D}; T^\ast \mathbb{D} \otimes \g \oplus \W \oplus (\Delta \otimes \V)_+ \right), \quad j = 1, 2,\] are two solutions to \eqref{eq:perturbed YM prelim}-\eqref{eq:perturbed D- convt}.   Assume also that they are gauge equivalent near $\p^-\mathbb{D}$, i.e. there exists $\mathbf{U} \in G^0 (\mathbb{D}, p)$ such that 
\[  \left(V^{(1)}, \Psi^{(1)}, \phi^{(1)}\right) = \left(V^{(2)}, \Psi^{(2)}, \phi^{(2)}\right) \cdot \mathbf{U} \quad \mbox{near $\p^- \mathbb{D}$}. \] 
If $(\mathcal{K}, \mathcal{J}, \mathcal{F})$ vanishes in $\mathbb{D} \setminus \mho$, $\mathbf{U} = \id$ in $\mho$ near $\p^- \mathbb{D}$, and $V^{(j)}$ is in the temporal gauge, then $\mathbf{U}$ is time-independent and $\left(V^{(1)}, \Psi^{(1)}, \phi^{(1)}\right)$ and $\left(V^{(2)}, \Psi^{(2)}, \phi^{(2)}\right)$ are gauge equivalent in $\mathbb{D}$.
 \end{proposition}

\subsection{Source-to-solution map} \label{subsec : StS}

Let $(\psi, A, \Phi)$ be a smooth critical point of $\mathscr{A}_{\mathrm{SM}}[\psi, A, \Phi]$.
In view of Proposition \ref{prop : solution in the relative gauge}
we may define the source-to-solution map 
\begin{align}\label{def_L}
\mathbf{L}_{(\psi, A, \Phi)} :      (\mathcal{K}, \mathcal{J}_1, \mathcal{J}_2, \mathcal{J}_3, \mathcal{F}) \longmapsto \mathscr{T}(\varphi, W,  \Upsilon) |_\mho,
\end{align}
where $(\varphi, W,  \Upsilon)$, with some $\mathcal J_0$, is the solution of \eqref{def_K_tilde}-\eqref{eqn : compatibility bg pt} 
for a small source \eqref{def_source_direct}.
Arguing as in Section 3.2 of \cite{CLOP3}, the key uniqueness result being now Proposition \ref{prop : uniqueness temporal},
we see that $\mathbf{L}_{(\psi, A, \Phi)}$ is uniquely determined by the data set $\mathcal{D}_{(\psi, A, \Phi)}(\mho)$ of the critical point.   
Theorem \ref{thm:main thm} reduces then to the following coefficient determination result for a system of nonlinear PDEs. 

\begin{theorem} \label{thm : StS to fields}
 	  	 Suppose the non-degeneracy conditions in Theorem \ref{thm:main thm} hold. The source-to-solution map $\mathbf{L}_{(\psi, A, \Phi)}$ determines the critical point $(\psi, A, \Phi)$ up to the gauge.
 	  \end{theorem}

\section{Reading principal symbols of linearized waves} \label{sec : linearization}

In view of the above reduction to Theorem \ref{thm : StS to fields}, we need to study the nonlinear system of wave equations \eqref{def_K_tilde}-\eqref{eqn : compatibility bg pt}.
We will use the technique of multiple linearizations originating from \cite{KLU}. In our case, three-fold linearization is used in the sense that we let the source \eqref{def_source_direct} depend on three independent parameters $\epsilon_1$, $\epsilon_2$ and $\epsilon_3$, and apply the cross derivative $\p_{\epsilon_1}\p_{\epsilon_2}\p_{\epsilon_3}$ to the source-to-solution map at such a source. Further, choosing the source to have a suitable structure as a conormal distribution, the cross derivative is also a conormal distribution, and we will compute its principal symbol (Proposition \ref{prop : principal symbol of 3-fold wave}).


\subsection{Linear dynamics}\label{subsec : linear dyn}

Let $(\psi, A, \Phi)$ be a smooth critical point of $\mathscr{A}_{\mathrm{SM}}[\psi, A, \Phi]$, and let $(\varphi, W,  \Upsilon, \mathcal J_0)$ solve \eqref{def_K_tilde}-\eqref{eqn : compatibility bg pt} with a source of the form $\left(\epsilon \mathcal K, \epsilon \mathcal J_1, \epsilon \mathcal J_2, \epsilon \mathcal J_3, \epsilon \mathcal F\right)$ for small $\epsilon > 0$. 
Differentiating \eqref{eq:perturbed D+ rL}-\eqref{eq:perturbed H rL} in $\epsilon$ at $\epsilon = 0$ yields a linear system  in 
	\begin{align}
(\varphi', W', \Upsilon') = \p_\epsilon (\varphi, W, \Upsilon)|_{\epsilon = 0}.
	\end{align}
Omitting the zeroth order terms, that is, the terms that do not contain any derivatives of $(\varphi', W', \Upsilon')$, the system reads
 \begin{align}	  
 	\Box_{A, \varrho} \varphi_{L}' + \frac{1}{2} d W' \bullet \psi_L      - \J_{\mathrm{YH}, L}(d \Upsilon'(e_a), \imath \Gamma^a \psi_R)  &= \slashed{D}_A\mathcal{K}_{L} 
 	\label{eq:SM pD+ linearized one-fold}\\
 	\Box_{A, \varrho} \varphi_{R}'   + \frac{1}{2} d W'   \bullet \psi_R      - \J_{\mathrm{YH}, R}(\imath \Gamma^a \psi_L, d \Upsilon'  (e_a))    &= \slashed{D}_A\mathcal{K}_{R}  
 	\label{eq:SM pD- linearized one-fold}\\
	\Box_{A, \Ad}W'     +    \J_{\mathrm{YMD}}(d\Upsilon',\Phi)       &=\mathcal{J}
	\label{eq:SM pYM linearized one-fold} \\	\label{eq:SM pH linearized one-fold}
	\Box_{A,\rho} \Upsilon'            &=\mathcal{F},
\end{align}  
where $\mathcal J = (\p_\epsilon \mathcal J_0|_{\epsilon=0},\mathcal J_1, \mathcal J_2, \mathcal J_3)$.
Here we used Lemma \ref{lemma:DAJ}.

All the three wave operators $\Box_{A, \cdot}$ have the same principal symbol 
    \begin{align*}
\sigma[\Box_{A, \cdot}](x,\xi) = \xi^\alpha \xi_\alpha,
    \end{align*} 
acting as a scalar multiplication. In other words, the system \eqref{eq:SM pD+ linearized one-fold}-\eqref{eq:SM pH linearized one-fold} consists of the scalar wave operator $\Box$ acting on each component in the leading order. 
It follows from \cite{Duistermaat-Hormander-FIO2} that if $(\varphi', W',  \Upsilon')$
is a conormal distribution in some region where the right-hand side of the system vanishes, 
then the principal symbol
    \begin{align*}
\sigma[(\varphi', W',  \Upsilon')] := (\varsigma, w, \upsilon)
    \end{align*}
satisfies the transport equation
\begin{align}
	\dot{\varsigma}_{L}+\varrho_{L*}(A_{\gamma}(\dot{\gamma}))\varsigma_{L}+\frac{1}{2}\varrho_{L*}(w)\cdot\dot{\gamma}\cdot\psi_{L}(\gamma)+\frac{1}{2}\J_{\mathrm{YH},L}(\upsilon,\dot{\gamma}\cdot\psi_{R}(\gamma))&=0\label{eq:ymhdyD+}\\
	\dot{\varsigma}_{R}+\varrho_{R*}(A_{\gamma}(\dot{\gamma}))\varsigma_{R}+\frac{1}{2}\varrho_{R*}(w)\cdot\dot{\gamma}\cdot\psi_{R}(\gamma)+\frac{1}{2}\J_{\mathrm{YH},R}(\upsilon,\dot{\gamma}\cdot\psi_{L}(\gamma))&=0 \label{eq:ymhdyD-}
	\\
\dot{w}_{\beta}+[A_{\gamma}(\dot{\gamma}),w_{\beta}]-\frac{1}{2}\dot{\gamma}_{\beta}\J_{\mathrm{YMH}}(\upsilon,\Phi(\gamma))&=0\label{eq:ymhdyYM}\\
\dot{\upsilon}+\rho_{*}(A_{\gamma}(\dot{\gamma}))\upsilon&=0.\label{eq:ymhdyH}
\end{align}
along lightlike geodesics $\gamma$. More precisely, as $(\varphi', W',  \Upsilon')$ is a conormal distribution satisfying a system of wave equations, its wave front set is invariant under the bicharacteristic flow, outside the support of the source. When projected on $\mathbb M$, the bicharacteristics are lightlike geodesics, that is, lines whose tangent vectors $v$ satisfy $v^\alpha v_\alpha = 0$. The transport equations hold along such lines. Adding zeroth order terms to \eqref{eq:SM pD+ linearized one-fold}-\eqref{eq:SM pH linearized one-fold} does not change the transport equations \eqref{eq:ymhdyD+}-\eqref{eq:ymhdyH}, so our omission of these is immaterial.
These facts are reviewed in detail in our previous work \cite{CLOP1}. 

On the level of the transport equations for principal symbols, the Yang--Mills and Higgs channels \eqref{eq:ymhdyYM}-\eqref{eq:ymhdyH} are decoupled from the Dirac channels \eqref{eq:ymhdyD+}-\eqref{eq:ymhdyD-}. This turns out to be important, since it will allow us to ignore the Dirac channels when recovering the Yang--Mills $A$ and Higgs $\Phi$ fields. 
Furthermore, we shall see that recovering the spinor field $\psi = \psi_L + \psi_R$
  relies solely on the terms
	\begin{align*}
\frac{1}{2}\varrho_{L*}(w)\cdot\dot{\gamma}\cdot\psi_{L}(\gamma), \quad \frac{1}{2}\varrho_{R*}(w)\cdot\dot{\gamma}\cdot\psi_{R}(\gamma). 
	\end{align*}
In fact, our strategy is to choose vanishing initial values for $\varsigma$
and $\upsilon$ whereas $w$ is in $Z(\g)$ initially. Then $\upsilon$ vanishes everywhere
and the system \eqref{eq:ymhdyD+}-\eqref{eq:ymhdyH} reduces to 
	\begin{align}
\label{transport_red} 
\dot{\varsigma} 
+ \varrho_{*}(A_{\gamma}(\dot{\gamma}))\varsigma
+ \frac{1}{2}\varrho_{*}(w)\cdot\dot{\gamma}\cdot\psi(\gamma)
&=0,
	\end{align} 
with $w$ coinciding with its initial value. 
Here it is not necessary to keep track of the splitting between the left and right Dirac channels. While the analysis of the transport equations becomes fairly simple with these reductions, the study of the interactions arising from the three-fold differentiation remains challenging. This is done in Section~\ref{sec:recovery}.  

Let us write the solution to \eqref{transport_red} in terms of parallel transports. 
First, we introduce the parallel transport $\mathbf{U}_\gamma^A$ on the principal bundle $\M \times G$ along a curve $\gamma: [0, T] \rightarrow M$. Namely, $\mathbf{U}_\gamma^A(t) = u(t)$, $t \in [0, T]$, is given by the solution to the following ODE
 \begin{align*}
 	\dot{u} + A_{\gamma}(\dot{\gamma}) \, u &= 0\\
 	u(0) &= \id.
 \end{align*} 
Then, given a vector space $\mathbb V$ and a linear representation $\kappa:G\to \GL(\mathbb{V})$, we obtain a parallel transport acting on $\mathbb{V}$ by setting 
\begin{align*}
\mathbf{P}_\gamma^{A, \kappa}(t) := \kappa(\mathbf{U}_\gamma^A(t)).
\end{align*}

Now the solution to \eqref{transport_red} takes the form
	\begin{align*} 
\varsigma(t)
&=
- \frac12  \P_\gamma^{A,\varrho}(t)\int_{0}^{t}  \varrho_{*}(w)\cdot\dot{\gamma} \cdot (\P_\gamma^{A,\varrho}(s))^{-1}\psi(\gamma(s))\,ds.
\end{align*} 
To simplify this formula further, let $w = b \otimes \omega \in Z(\g) \otimes \R^4$, 
and write 
	\begin{align}\label{truncated_ray_trans}
I_{\gamma}(t) 
= 
\P_\gamma^{A,\varrho}(t) \int_0^t (\P_\gamma^{A,\varrho}(s))^{-1} \psi(\gamma(s))\,ds.
	\end{align}
Then we rewrite $\varsigma(t)$ as 
    \begin{align}\label{transport_sol}
\varsigma(t)
=
- \frac12
\omega \cdot \dot \gamma \cdot \varrho_{*}(b) I_\gamma(t).
    \end{align}
We will show that for certain vectors $v$ the quantity $v \cdot \varsigma(t)$ can be recovered
via the three-fold linearization scheme.

\subsection{Three-fold linearization} \label{subsec : linearization}

The strategy outlined Section \ref{subsec : linear dyn} uses sources only in the Yang--Mills and Higgs channels, not in the Dirac channels. For this reason, we will take $\mathcal K_L = 0$ and $\mathcal K_R = 0$ in what follows. Then the constraint \eqref{def_K_tilde} is satisfied. 

Let $(\psi, A, \Phi)$ be a smooth critical point of $\mathscr{A}_{\mathrm{SM}}[\psi, A, \Phi]$, and let $(\varphi, W,  \Upsilon, \mathcal J_0)$ solve \eqref{eq:perturbed D+ rL}-\eqref{eqn : compatibility bg pt} with a source of the form 
	\begin{align}\label{def_vecJF}
\left(\mathcal{K}, \mathcal{J}_1, \mathcal{J}_2, \mathcal{J}_3, \mathcal{F}\right) = \left(
0, \epsilon \cdot	\vec{\mathcal{J}}_1, 
\epsilon \cdot	\vec{\mathcal{J}}_2, 
 \epsilon \cdot	\vec{\mathcal{J}}_3,
\epsilon \cdot \vec{\mathcal{F}}\right),
	\end{align}
where 
\begin{align}\label{eq : vector sources}\left\{
	\begin{aligned}\epsilon &:= \left(\epsilon_{(1)}, \epsilon_{(2)}, \epsilon_{(3)}\right) \in \R^3 \\
 	\vec{\mathcal{J}}_k &:= \left(\mathcal{J}_{(1),k}, \mathcal{J}_{(2),k}, \mathcal{J}_{(3),k}\right) \in \left(C_c^\infty(\mho,\g)\right)^3, \quad k=1,2,3, \\
    \vec{\mathcal{F}} &:= \left(\mathcal{F}_{(1)}, \mathcal{F}_{(2)}, \mathcal{F}_{(3)}\right) \in \left(C_c^\infty(\mho,  \W)\right)^3.
\end{aligned}\right.
 \end{align}
Analogously to our previous work \cite{CLOP3},
the linear and the quadratic terms with no derivatives hitting $(\varphi, W,  \Upsilon)$
play no role in our subsequent computations. As the formulas for two and three-fold linearizations are fairly complicated, we will omit writing these terms. Without them, equations \eqref{eq:perturbed D+ rL}-\eqref{eq:perturbed H rL} read
 \begin{align}
 	\Box_{A, \varrho} \varphi_L + \frac{1}{2} d W \bullet \psi_L - \J_{\mathrm{YH}, L}(d \Upsilon(e_a), \imath \Gamma^a \psi_R)    &=  	  N^{\mathrm{DL}}
 	\label{eq:SM pD+ principal}\\
 	\Box_{A, \varrho} \varphi_R + \frac{1}{2} d W \bullet \psi_R - \J_{\mathrm{YH}, R}(\imath \Gamma^a \psi_L, d \Upsilon(e_a))  &=     N^{\mathrm{DR}}
 	\label{eq:SM pD- principal}
 	\\
 	\Box_{A, \Ad}W + \J_{\mathrm{YMH}}(d\Upsilon, \Phi)        &=\mathcal{J}+   N^{\mathrm{YM}}
\label{eq:SM pYM principal} \\	
\Box_{A,\rho} \Upsilon      &=\mathcal{F}    +  	N^{\mathrm{H}} ,
\label{eq:SM pH principal}  
 \end{align}  where  
\begin{align*}  
-N^{\mathrm{DL}} &:=\frac{1}{2} d W \bullet  \varphi_L  + 2 \star (\varrho_\ast(W) \wedge \star d \varphi_L)  \\
    &\quad +  \star (\varrho_\ast(W) \wedge \star (\varrho_\ast(W) \varphi_L)) + \frac{1}{4} [W, W] \bullet \varphi_L     -  \slashed{D}_W \J_{\mathrm{YH}, L}(\Upsilon, \varphi_R)
\end{align*} 
$N^{\mathrm{DR}}$ is defined analogously, swapping $L$ and $R$, and
	\begin{align*}	-N^{\mathrm{YM}} &:= \frac{1}{2} d^\ast [W, W] + \star [W, \star d W] + \frac{1}{2} \star [W, \star [W, W]]   
	\\
	&\quad + \J_{\mathrm{YMH}}(d \Upsilon,\Upsilon) + \J_{\mathrm{YMH}}(\rho_{*}(W)\Upsilon,\Upsilon)
	\\ 
-N^{\mathrm{H}} &:= 2 \star (\rho_\ast (W) \wedge \star d \Upsilon) + \star(\rho_{*}(W)\wedge\star \rho_{*}(W)\Upsilon)  	+ |\Upsilon|^2_{\W} \Upsilon. 
\end{align*}

We now perform the third order linearization scheme on \eqref{eq:SM pD+ principal}-\eqref{eq:SM pH principal}. That is, we differentiate the system in $\epsilon_1$, $\epsilon_2$ and $\epsilon_3$ at $0$. We write
for $X = \varphi_L, \varphi_R, W, \Upsilon$,
\begin{align}\label{eq: linearized fields}\left\{\begin{aligned}
 X_{(j)}	&:=	\p_{\epsilon_{(j)}}   X |_{\epsilon  = 0}  \\
 X_{(jk)}	&:=  \p_{\epsilon_{(j)}}\p_{\epsilon_{(k)}}   X  |_{\epsilon = 0}  \\	
 X_{(123)}	&:=  \p_{\epsilon_{(1)}}\p_{\epsilon_{(2)}}\p_{\epsilon_{(3)}}   X   |_{\epsilon = 0} .
\end{aligned}\right.
\end{align} 
Then $(\varphi_{(j)}, W_{(j)}, \Upsilon_{(j)})$ satisfies \eqref{eq:SM pD+ linearized one-fold}-\eqref{eq:SM pH linearized one-fold} with $(\mathcal{J}, \mathcal{F}) = (\mathcal{J}_{(j)}, \mathcal{F}_{(j)})$
and $(\mathcal K_L, \mathcal K_R) = 0$.
The two-fold equations read 
 \begin{align}
 	\Box_{A, \varrho} \varphi_{L, (jk)} + \frac{1}{2} d W_{(jk)} \bullet \psi_{L}  - \J_{\mathrm{YH}, L}(d \Upsilon_{(jk)}(e_a), \imath \Gamma^a \psi_R)  	  &= N^{\mathrm{DL}}_{(jk)}
 	\label{eq:SM pD+ linearized 2-fold}\\
 	\Box_{A, \varrho} \varphi_{R, (jk)} + \frac{1}{2} d W_{(jk)} \bullet \psi_{R} - \J_{\mathrm{YH}, R}(\imath \Gamma^a \psi_L, d \Upsilon_{(jk)}(e_a))   &= N^{\mathrm{DR}}_{(jk)} 
 	\label{eq:SM pD- linearized 2-fold}
 	\\
	\Box_{A, \Ad}W_{(jk)}   +   \J_{\mathrm{YMH}}(d\Upsilon_{(jk)},\Phi)       &=N^{\mathrm{YM}}_{(jk)}
	\label{eq:SM pYM linearized 2-fold} \\	
	\Box_{A,\rho} \Upsilon_{(jk)}           &=N^{\mathrm{H}}_{(jk)},
	\label{eq:SM pH linearized 2-fold}  
\end{align}  where 
\begin{align*}  
	-N^{\mathrm{DL}}_{(jk)} &:= 2\star \left(\varrho_\ast(W_{(j)}) \wedge \star d \varphi_{L, (k)}\right)    + \frac{1}{2} d W_{(j)} \bullet \varphi_{L, (k)} 
	  -  \J_{\mathrm{YH}, L}(d \Upsilon_{(j)}(e_a), \imath \Gamma^a \varphi_{R, (k)}) \\
	&\quad  + 2\star \left(\varrho_\ast(W_{(k)}) \wedge \star d \varphi_{L, (j)}\right)    + \frac{1}{2} d W_{(k)} \bullet \varphi_{L, (j)}  
	  - 	\J_{\mathrm{YH}, L}(d \Upsilon_{(k)}(e_a), \imath \Gamma^a \varphi_{R, (j)})    
\end{align*} 
and $N^{\mathrm{DR}}_{(jk)}$ is analogous. Moreover, 
\begin{align*} 
		-N^{\mathrm{YM}}_{(jk)} &:= \frac{1}{2} d^\ast [W_{(j)}, W_{(k)}] + \star [W_{(j)}, \star d W_{(k)}]   + \J_{\text{YMH}}(d \Upsilon_{(j)},\Upsilon_{(k)}) \\
	&\quad + \frac{1}{2} d^\ast [W_{(k)}, W_{(j)}] + \star [W_{(k)}, \star d W_{(j)}]   + \J_{\text{YMH}}(d \Upsilon_{(k)},\Upsilon_{(j)})
	\\ 
	-N^{\mathrm{H}}_{(jk)} &:= 2 \star (\rho_\ast (W_{(j)}) \wedge \star d \Upsilon_{(k)}) + 2 \star (\rho_\ast (W_{(k)}) \wedge \star d \Upsilon_{(j)}).
\end{align*}

Finally, the three-fold equations take the form 
 \begin{align}
 	\Box_{A, \varrho} \varphi_{L, (123)} + \frac{1}{2} d W_{(123)} \bullet \psi_L  - \J_{\mathrm{YH}, L}(d \Upsilon_{(123)}(e_a), \imath \Gamma^a \psi_R)   &= N^{\mathrm{DL}}_{(123)}
 	\label{eq:SM pD+ linearized 3-fold}\\
 	\Box_{A, \varrho} \varphi_{R, (123)} + \frac{1}{2} d W_{(123)} \bullet \psi_R  - \J_{\mathrm{YH}, R}(\imath \Gamma^a \psi_L, d \Upsilon_{(123)}(e_a))   &= N^{\mathrm{DR}}_{(123)}  
 	\label{eq:SM pD- linearized 3-fold} 
 	\\
	\Box_{A, \Ad}W_{(123)} + \J_{\mathrm{YMH}}(d\Upsilon_{(123)},\Phi)        &=N^{\mathrm{YM}}_{(123)}
	\label{eq:SM pYM linearized 3-fold} \\	
	\Box_{A,\rho} \Upsilon_{(123)}           &=N^{\mathrm{H}}_{(123)},
	\label{eq:SM pH linearized 3-fold}  
\end{align}  where 
\begin{align*}  
-N^{\mathrm{DL}}_{(123)} &:= \frac{1}{2}\sum_\pi \bigg( 2\star \left(\varrho_\ast(W_{(\pi(1)\pi(2))}) \wedge \star d \varphi_{L, (\pi(3))}\right)    +  2\star \left(\varrho_\ast(W_{(\pi(1))}) \wedge \star d \varphi_{L, (\pi(2)\pi(3))}\right)  
\\ &\quad  +     \frac{1}{2} d W_{(\pi(1)\pi(2))} \bullet \varphi_{L, (\pi(3))} +     \frac{1}{2} d W_{(\pi(1))} \bullet \varphi_{L, (\pi(2)\pi(3))}  
	\\
  &\quad-	\J_{\mathrm{YH}, L}(d \Upsilon_{(\pi(1)\pi(2))}(e_a), \imath \Gamma^a \varphi_{R, (\pi(3))})  -	\J_{\mathrm{YH}, L}(d \Upsilon_{(\pi(1))}(e_a), \imath \Gamma^a \varphi_{R, (\pi(2)\pi(3))})   \\
  &\quad  - 2\J_{\mathrm{YH}, L}(\Upsilon_{(\pi(1))}, \J_{\mathrm{YH}, R}(\varphi_{L, (\pi(2))}, \Upsilon_{\pi(3)})) - 2 W_{(\pi(1))} \bullet \J_{\mathrm{YH}, L}(\Upsilon_{(\pi(2))}, \varphi_{R, (\pi(3))})  \\ &\quad   + 2\star \left(\varrho_\ast(W_{(\pi(1))}) \wedge \star (\varrho_\ast(W_{(\pi(2))}) \varphi_{(\pi(3))}))\right)  
  + \frac{1}{2} [W_{(\pi(1))}, W_{(\pi(2))}] \bullet \varphi_{(\pi(3))}  \bigg)
\end{align*} 
with $N^{\mathrm{DR}}_{(123)}$ analogous, and	\begin{align*}
  -N^{\mathrm{YM}}_{(123)} &:= \frac{1}{2} \sum_\pi \bigg(\frac{1}{2} d^\ast [W_{(\pi(1)\pi(2))}, W_{(\pi(3))}]   + \frac{1}{2} d^\ast [W_{(\pi(1))}, W_{(\pi(2)\pi(3))}] \\ 
  &\quad +\star [W_{(\pi(1)\pi(2))}, \star d W_{(\pi(3))}] +\star [W_{(\pi(1))}, \star d W_{(\pi(2)\pi(3))}] \\ 
  &\quad +\J_{\text{YMH}}(d \Upsilon_{(\pi(1)\pi(2))},\Upsilon_{(\pi(3))}) +\J_{\text{YMH}}(d \Upsilon_{(\pi(1))},\Upsilon_{(\pi(2)\pi(3))})\\ 
  &\quad +  2 \star [W_{(\pi(1))}, \star [W_{(\pi(2))}, W_{(\pi(3))}]]+2\J_{\text{YMH}}(\rho_{*}(W_{(\pi(1))})\Upsilon_{(\pi(2))},\Upsilon_{(\pi(3))})\bigg)
\\
  -N^{\mathrm{H}}_{(123)} &:= \sum_\pi \bigg( \star (\rho_\ast (W_{(\pi(1)\pi(2))}) \wedge \star d \Upsilon_{(\pi(3))})+  \star (\rho_\ast (W_{(\pi(1))}) \wedge \star d \Upsilon_{(\pi(2)\pi(3))})\\
  &\quad +\star(\rho_{*}(W_{(\pi(1))})\wedge\star \rho_{*}(W_{(\pi(2))})\Upsilon_{(\pi(3))})  	+ \langle \Upsilon_{(\pi(1))}, \Upsilon_{(\pi(2))} \rangle \Upsilon_{(\pi(3))}    \bigg).
\end{align*} 
 
\subsection{Special solutions with conormal structure} \label{subsec : from nonlinear to multilinear}

Following our previous work in the series \cite{CLOP1, CLOP2, CLOP3}, we shall make use of a class of solutions to \eqref{eq:SM pD+ principal}-\eqref{eq:SM pH principal} that are conormal distributions with suitable geometric properties. In order to describe the geometric setup we introduce some notation. 

We set $v^\flat = v_\alpha dx^\alpha$ for a vector $v = v^\alpha \p_{x^\alpha}$,
and $\xi^\sharp = \xi^\alpha \p_{x^\alpha}$ for a covector $\xi = \xi_\alpha dx^\alpha$.
Consider the sets $\mathbb D$ and $\mho$ defined by \eqref{def_diamond} and \eqref{def_mho}, see also Figure \ref{fig_D}. For $y \in \mathbb D$ there is $x \in \mho$ such that $y - x$ is a lightlike future pointing vector.  
After a rotation, we may assume that $x = y + \ell \xi_{(1)}^\sharp$ where $\ell > 0$ and 
	\begin{align*}
\xi_{(1)} = (1, 1, 0, 0). 
	\end{align*}
We write
	\begin{align*}
\eta = (1, -a(r), r, 0),
	\end{align*}
where $a(r) = \sqrt{1 - r^2}$ and $r \ne 0$ satisfies $|r| \le 1$. The point $z = y - \ell \eta^\sharp$ is in $\mho$ when $|r|$ is small.

We choose 
	\begin{align*}
\xi_{(2)} = (1, a(s), s, 0), \quad \xi_{(3)} = (1, a(s), -s, 0),
	\end{align*}
with $s > 0$ small enough so that $x_{(i)} \in \mho$ where $x_{(i)} = y + \ell \xi_{(i)}^\sharp$, $i = 2,3$. We write also $x_{(1)} = x$. The lines 
	\begin{align}\label{def_gamma_i}
\gamma_{(i)}(t) = y + (\ell - t) \xi_{(i)}^\sharp, \quad i=1,2,3,
	\end{align}
satisfy $\gamma_{(i)}(0) = x_{(i)}$ and $\gamma_{(i)}(\ell) = y$,
and the line
	\begin{align}\label{def_gamma_4}
\gamma_{(4)}(t) = y - t\eta^\sharp,
	\end{align}
satisfies $\gamma_{(4)}(0) = y$ and $\gamma_{(4)}(\ell) = z$.

We choose the components of $\vec{\mathcal{J}}_k$ and $\vec{\mathcal{F}}$, see \eqref{eq : vector sources}, as follows
\begin{align} \label{eqn : choice of source} 
	\left\{
 \begin{aligned}	
     \mathcal{J}_{(i), k}(s) &:= b_{(i)}\omega_{(i), k} \chi_{(i)} \delta_{x_{(i)}}, \quad k=1,2,3,    
     \\ 
     \mathcal{F}_{(i)}(s) &:= \upsilon_{(i)} \chi_{(i)} \delta_{x_{(i)}}  ,  
 \end{aligned} 
     \right.
\end{align}
 where 
 \begin{itemize}
 	\item $\omega_{(i),k} \in \R$, $b_{(i)} \in  \mathfrak{g}$, $\upsilon_{(i)} \in \mathcal{W}$;
 	\item $\delta_{x_{(i)}}$ is the Dirac delta distribution at $x_{(i)}$;
 	
 	\item $\chi_{(i)}$ is a microlocalization near the two points $(x_{(i)}, \pm \xi_{(i)})$ $i = 1, 2, 3$; 
 	\item the principal symbol $\sigma[\chi_{(i)}]$ is positively homogeneous of degree $q \leq -9$ (the homogeneity guarantees that the sources $\mathcal{J}_{(i), k}$ and $\mathcal{F}_{(i)}$ satisfy the smoothness requirement \eqref{def_source_direct});
 	\item $\mho_{(i)} \cap \mathcal J^+(\mho_{(l)}) = \emptyset$ for all $i \ne l$, where $\mho_{(i)} \subset \mho$ is a neighbourhood  of $x_{(i)}$, covering $$\supp \mathcal{J}_{(i), 1}  \cup \supp \mathcal{J}_{(i), 2}  \cup \supp \mathcal{J}_{(i), 3} \cup \supp   \mathcal{F}_{(i)}   ,$$ and $\mathcal J^+(\mho_{(l)})$ is the causal future of $\mho_{(l)}$;

 	\item $\hat \mho_{(i)} \cap \Gamma_{(l)} = \emptyset$ for all $i \ne l$, where
 	\begin{align*}
 		\hat \mho_{(i)} &:= \{(t,x') \in \R^{1+3}: (\tilde t, x') \in \mho_{(i)} \text{ for some $\tilde t \le t$} \},
 		\\
 		\Gamma_{(i)} &:= \{x_{(i)} + t\xi^\sharp : t \in \R,\ (x_{(i)}, \xi) \in \mathrm{WF}(\chi_{(i)}) \}.
 	\end{align*}   	
 \end{itemize}

The following is an analogue of \cite[Proposition 5]{CLOP3}, with essentially the same proof.

\begin{proposition}\label{prop : principal symbol of 3-fold wave}
Consider a source of the form \eqref{def_vecJF} where $\mathcal J_{(i),k}$ and $\mathcal F_{(i)}$ are as in \eqref{eqn : choice of source}.
Then 
\begin{align}
\label{eqn : principal symbol of 3-fold wave}	\sigma\left[ \left.\p^3_{\epsilon_{(1)}\epsilon_{(2)}\epsilon_{(3)}} \mathbf{L}_{(\psi, A, \Phi)}(0, \mathcal{J}_1, \mathcal{J}_2, \mathcal{J}_3, \mathcal{F}) \right|_{\epsilon = 0} \right] (z, \eta) 
= (\varsigma, w, \upsilon),
\end{align} 
where, letting $(\varphi_{(123)}, W_{(123)}, \Upsilon_{(123)})$ to be the solution of \eqref{eq:SM pD+ linearized 3-fold}-\eqref{eq:SM pH linearized 3-fold}, 
\begin{align*}
	\varsigma &= \sigma[\varphi_{(123)}] (z, \eta) \\ 
		w_{\beta} &= - \frac{\eta_\beta}{\eta_0} \sigma[W_{(123), 0}](z, \eta) + \sigma[W_{(123), \beta}](z, \eta)\\
		\upsilon &= \sigma[\Upsilon_{(123)}](z, \eta).
\end{align*} 
\end{proposition}

The form of the components $w_\beta$ is due to the temporal gauge map $\mathscr T$ being a factor of the source-to-solution map \eqref{def_L}. 
Varying $\omega_{(i),k} \in \R$, $b_{(i)} \in  \mathfrak{g}$ and $\upsilon_{(i)} \in \mathcal{W}$, as well as $y \in \mathbb D \setminus \mho$ and the directions $\eta$ and $\xi_{(i)}$, we will recover from \eqref{eqn : principal symbol of 3-fold wave} the critical point $(\psi, A, \Phi)$, up to the gauge.

\section{Recovery of the critical point} \label{sec:recovery}

\subsection{Yang--Mills and Higgs components}

\label{subsec : YM recovery}
We begin by reconstructing the Yang--Mills and Higgs components $(A, \Phi)$ of the critical point $(\psi, A, \Phi)$ from the knowledge of \eqref{def_L}. This is achieved by a straightforward reduction to the  inverse problem for the Yang--Mills--Higgs system which we solved in \cite{CLOP3}. 

The perturbed Yang--Mills--Higgs system takes the form
\begin{align}  
	D_{A}^*D_{A}W  + \mathcal{N}_{A}(W)+\sum_{i=1}^{3}\mathcal{M}_{A,\Phi}^{i}(W,\Upsilon)    &= \mathcal{J}
	\label{eq:perturbed YM in YMH} \\	
	\Box_{A,\rho} \Upsilon+ \sum_{i=1}^{3} 	\mathcal{O}_{A, \Phi}^i(W, \Upsilon)   &=\mathcal{F}.
	\label{eq:perturbed H in YMH} 
\end{align} where $\mathcal{M}_{A,\Phi}^{i}$ and $\mathcal{O}_{A, \Phi}^i(W, \Upsilon)$ are the same as in \eqref{eq:perturbed YM prelim}-\eqref{eq:perturbed H prelim}.
In \cite{CLOP3} we associated to \eqref{eq:perturbed YM in YMH}-\eqref{eq:perturbed H in YMH} a source-to-solution map $$ \mathbf{L}_{(A, \Phi)} : \left(\mathcal{J}_1, \mathcal{J}_2, \mathcal{J}_3, \mathcal{F}\right) \longmapsto \mathscr{T}(W, \Upsilon)|_\mho,$$ 
where $W, \Upsilon$ solves \eqref{eq:perturbed YM in YMH}-\eqref{eq:perturbed H in YMH} in the relative Lorenz gauge. 

Writing $\tilde{\mathbf{L}}_{(\psi, A, \Phi)}$ 
for the restriction of $\mathbf{L}_{(\psi, A, \Phi)}$ to the Yang--Mills and Higgs channels, we have using the notation in \eqref{eqn : principal symbol of 3-fold wave}, 
\begin{multline*}\sigma\left[ \left.\p^3_{\epsilon_{(1)}\epsilon_{(2)}\epsilon_{(3)}} \tilde{\mathbf{L}}_{(\psi, A, \Phi)}(0, \mathcal{J}_1, \mathcal{J}_2, \mathcal{J}_3, \mathcal{F}) \right|_{\epsilon = 0} \right]  
= (w, \upsilon)
\\
= \sigma\left[ \left.\p^3_{\epsilon_{(1)}\epsilon_{(2)}\epsilon_{(3)}} \mathbf{L}_{(A, \Phi)}(\mathcal{J}_1, \mathcal{J}_2, \mathcal{J}_3, \mathcal{F}) \right|_{\epsilon = 0} \right],\end{multline*}
where the source $(\mathcal{J}_1, \mathcal{J}_2, \mathcal{J}_3, \mathcal{F})$ is as in Proposition \ref{prop : principal symbol of 3-fold wave}.
This is because the Yang--Mills and Higgs channels are decoupled from the Dirac channels on the level of the transport equations \eqref{eq:ymhdyD+}-\eqref{eq:ymhdyH}, and because the leading order nonlinear interactions $N^{\mathrm{YM}}_{(jk)}$, $N^{\mathrm{H}}_{(jk)}$, $N^{\mathrm{YM}}_{(123)}$ and $N^{\mathrm{H}}_{(123)}$ do not depend on the linearized Dirac fields $\varphi_{(j)}$ and $\varphi_{(jk)}$.

In \cite{CLOP3} we proved  that the symbols
\[\sigma\left[ \left.\p^3_{\epsilon_{(1)}\epsilon_{(2)}\epsilon_{(3)}} \mathbf{L}_{(A, \Phi)}(\mathcal{J}_1, \mathcal{J}_2, \mathcal{J}_3, \mathcal{F}) \right|_{\epsilon = 0} \right],\] 
for sources $(\mathcal{J}_1, \mathcal{J}_2, \mathcal{J}_3, \mathcal{F})$ as in Proposition \ref{prop : principal symbol of 3-fold wave},
determine the orbit  
	\begin{align*}
(A, \Phi) \cdot G^0(\mathbb D, p) 
= 
\{(A, \Phi) \cdot \mathbf{U} \mid \mathbf{U} \in G^0(\mathbb D, p)\},
	\end{align*}
the optimal result in view of the gauge invariance.
Let us choose a point $(\tilde A, \tilde \Phi)$ on the orbit. Then there is $\tilde \psi$ such that $(\tilde \psi, \tilde A, \tilde \Phi) \sim (\psi, A, \Phi)$. This gauge equivalence implies
	\begin{align*}
\mathcal{D}_{(\tilde \psi, \tilde A, \tilde \Phi)}(\mho) = \mathcal{D}_{(\psi, A, \Phi)}(\mho),
	\end{align*}
and, to simplify the notation, we replace
$(\psi, A, \Phi)$ by $(\tilde \psi, \tilde A, \tilde \Phi)$.
We will show below that knowing the symbols \eqref{eqn : principal symbol of 3-fold wave} and $(A,\Phi)$ determines $\psi$, not just up to a gauge but fully. However, due to the change of notation, this means that we have then recovered $(\psi, A, \Phi)$ up to a gauge. 

\subsection{Interactions in the Dirac channel}\label{subsec : EM signal in linearized system}

Let us begin by explaining the rationale behind the further structure that we impose on the source. 
Recall that if the initial conditions for the system of transport equations \eqref{eq:ymhdyD+}-\eqref{eq:ymhdyH} are chosen so that $\varsigma(0) = 0$, $\upsilon(0) = 0$ and $w(0) \in Z(\g)$, then $\upsilon = 0$ and $w = w(0)$ everywhere, and the system reduces to the single equation \eqref{transport_red}.
The leading order nonlinear interactions $N^{\mathrm{YM}}_{(jk)}$, $N^{\mathrm{H}}_{(jk)}$, $N^{\mathrm{YM}}_{(123)}$ and $N^{\mathrm{H}}_{(123)}$ in the Yang--Mills and Higgs channels vanish, on the level of the principal symbols, since all the terms in these expressions contain $\Upsilon$ factors or a commutator. The former vanish since $\upsilon$ = 0 and latter since $w \in Z(\g)$. Thus we only need to analyze the interactions $N^{\mathrm{DL}}_{(jk)}$, $N^{\mathrm{DR}}_{(jk)}$, $N^{\mathrm{DL}}_{(123)}$ and $N^{\mathrm{DR}}_{(123)}$ in the Dirac channel. 

One may wonder if taking $w(0) = 0$ and nonzero $\varsigma(0)$ would lead to a simpler analysis. In this case, $w = 0$ everywhere and the leading order nonlinear interactions vanish in all channels. Indeed, the interactions in the Yang--Mills and Higgs channels vanish as above, and the interactions in the Dirac channel vanish since all the terms in $N^{\mathrm{DL}}_{(jk)}$, $N^{\mathrm{DR}}_{(jk)}$, $N^{\mathrm{DL}}_{(123)}$ and $N^{\mathrm{DR}}_{(123)}$ contain $W$ or $\Upsilon$. Thus we would need to study the subleading interactions. We have opted not to do this.

We proceed to the detailed analysis of the interactions in the Dirac channel. Consider a source of the form \eqref{def_vecJF} where $\mathcal J_{(i),k}$ and $\mathcal F_{(i)}$ are as in \eqref{eqn : choice of source}.
Further, choose $b_{(i)}$, $\omega_{(i),k}$ and $\upsilon_{(i)}$ in \eqref{eqn : choice of source}
as follows 
\begin{equation} \label{eqn : special sources} \left\{ \begin{aligned} 
& b_{(1)} = b_{(2)} = b_{(3)} =: b \in Z(\g),
\\ & \omega_{(1),2} = \omega_{(2),2} = \omega_{(3),2} = 1, 
\\ & \omega_{(1),k} = \omega_{(2),k} = \omega_{(3),k} = 0,
\quad k = 1,3,
\\ &  \upsilon_{(1)}  = \upsilon_{(2)}  = \upsilon_{(3)}  = 0.\end{aligned}\right.\end{equation} 
At the level of principal symbols, the relative Lorenz gauge condition $D_A^* W = 0$ reduces to $\xi^\alpha \sigma[W_\alpha](x,\xi)=0$. 
As in our previous works \cite{CLOP2, CLOP3}, this implies that 
$\sigma[W_{(j)}](y,\xi_{(j)})]$ is proportional to $\omega_{(j)} b$ where
	\begin{align*}
\omega_{(1)} = dx^2, \quad \omega_{(2)} = s dx^0 + dx^2, \quad  \omega_{(3)} = - s dx^0 + dx^2.
	\end{align*}

Writing
	\begin{align}\label{def_kappa}
\kappa_{(1)}  := 1 - \frac{1 + a(r)}{1 - a(s)},
\quad
\kappa_{(2)}  := \frac{1+a(r)}{2(1-a(s))} + \frac{r}{2s},
\quad
\kappa_{(3)}  := \frac{1+a(r)}{2(1-a(s))} - \frac{r}{2s},
	\end{align}
and $\eta_{(i)} = \kappa_{(i)} \xi_{(i)}$, it  holds that
\[\eta = \eta_{(1)} + \eta_{(2)} + \eta_{(3)} = \kappa_{(1)} \xi_{(1)} + \kappa_{(2)} \xi_{(2)} + \kappa_{(3)} \xi_{(3)}.\]
As the microlocal cut-off $\chi_{(i)}$ in \eqref{eqn : choice of source} is homogeneous of degree $q$, we have
	\begin{align}
\label{eqn : 1-fold principal symbol}
\sigma[W_{(i)}](y, \eta_{(i)})  = \alpha_{(i)} |\kappa_{(i)}|^{q - 1} \omega_{(i)} b,
	\end{align}
where $\alpha_{(i)} \ne 0$ can be viewed as a volume factor, independent of the critical point $(\psi, A, \Phi)$.
The precise form of $\alpha_{(i)}$, given in \cite[(65)]{CLOP1},
is not important for our purposes. The principal symbol $\sigma[\varphi_{(i)}](y, \eta_{(i)})$ scales in the same way.

We write for $X = \varphi, W, \Upsilon$
\begin{align*}
	\hat{X}_{(j)} &:= \alpha_{(j)}^{-1} |\kappa_{(j)}|^{1-q} \sigma[X_{(j)}](y, \eta_{(j)}),
\end{align*}
and for $X = N^{\mathrm{D}}, N^{\mathrm{YM}}, N^{\mathrm{H}}$ and $Y = \varphi, W, \Upsilon$
	\begin{align*}
		\hat{X}_{(kl)} &:= (\alpha_{(k)} \alpha_{(l)})^{-1} |\kappa_{(k)}\kappa_{(l)}|^{1-q} \sigma[X_{(kl)}](y, \eta_{(kl)}), \\
			\hat{Y}_{(kl)} &:= - \imath (\alpha_{(k)} \alpha_{(l)})^{-1} |\kappa_{(k)}\kappa_{(l)}|^{1-q} \sigma[Y_{(kl)}](y, \eta_{(kl)}), 
	\end{align*}
where $\eta_{(kl)} = \eta_{(k)} + \eta_{(l)}$. Further, we write 
	\begin{align*}
\hat{N}^{\mathrm{D}}_{(123)} &:= (\alpha_{(1)}\alpha_{(2)}\alpha_{(3)})^{-1} |\kappa_{(1)}\kappa_{(2)}\kappa_{(3)}|^{1-q} \sigma[N^{\mathrm{D}}_{(123)}](y, \eta). 
	\end{align*}
As explained above, with our choice of the source, 
$\hat \Upsilon_{(j)} = 0$ and $\hat W_{(j)} = \omega_{(j)} b$ where $b \in Z(\g)$. Therefore $\hat N^{\mathrm{YM}}_{(jk)}$ and $\hat N^{\mathrm{H}}_{(jk)}$ vanish, and so do $\hat W_{(jk)}$ and $\hat \Upsilon_{(jk)}$. 
This again implies that $\hat{N}^{\mathrm{YM}}_{(123)}$ and $\hat{N}^{\mathrm{H}}_{(123)}$ vanish.
Due to $\hat \Upsilon_{(j)} = 0$ and $\hat \Upsilon_{(jk)} = 0$ there is no need to keep track of the splitting between the left and right Dirac channels. We have
\begin{align*} 
	  \hat N^{\mathrm{D}}_{(kl)} &:=	 - 2 \imath \star \left(\varrho_\ast(\hat W_{(k)}) \wedge \star \eta_{(l)}  \hat \varphi_{(l)}\right) - 2 \imath \star \left(\varrho_\ast(\hat W_{(l)}) \wedge \star \eta_{(k)} \hat \varphi_{(k)}\right)  \\&\quad  - \frac{\imath}{2} (\eta_{(k)} \wedge \hat W_{(k)}) \bullet \hat \varphi_{(l)}   - \frac{\imath}{2} (\eta_{(l)} \wedge \hat W_{(l)}) \bullet \hat \varphi_{(k)},
\end{align*}
and
	\begin{align}\label{ND_123_split}
\hat{N}^{\mathrm{D}}_{(123)} 
= 
\hat N^{\mathrm{D} (1+1+1)}_{(123)}
+ 
\hat N^{\mathrm{D}\wedge}_{(123)}
+ 
\hat N^{\mathrm{D}\bullet}_{(123)},
	\end{align}
where
\begin{align*} 	  	\hat N^{\mathrm{D} (1+1+1)}_{(123)} &:= -\sum_{\pi \in S_3}\star (\varrho_\ast(\hat W_{(\pi(1))}) \wedge \star (\varrho_\ast(\hat W_{(\pi(2))}) \hat \varphi_{(\pi(3))}))  \\
	\hat N^{\mathrm{D}\wedge}_{(123)} &:=	 - \imath \sum_{\pi \in S_3}\star \left(\varrho_\ast(\hat W_{(\pi(1))}) \wedge \star \eta_{(\pi(2)\pi(3))} \hat\varphi_{(\pi(2)\pi(3))}\right) 
	\\ 	\hat N^{\mathrm{D}\bullet}_{(123)} &:=  - \frac{\imath}{4} \sum_{\pi \in S_3}(\eta_{(\pi(1))} \wedge \hat W_{(\pi(1))}) \bullet \hat\varphi_{(\pi(2)\pi(3))}.
\end{align*} 

\subsubsection{Two-fold interactions}

We write $\breve{b}_{(i)} := \varrho_\ast\left(  b_{(i)}\right)$,
and consider the four terms in $\hat N^{\mathrm{D}}_{(kl)}$. While $\breve{b}_{(i)}$ are the same for $i=1,2,3$, we keep the subscript in the notation for the moment.
It acts as bookkeeping device since the indices are permuted systematically. 

Using 
	\begin{align*}
\star \left(  dx^\alpha \wedge \star dx^\beta \right) = - g^{\alpha\beta}.
	\end{align*}
the first term expands to
\begin{align*}
- 2 \imath \star \left(\varrho_\ast(\hat W_{(k)}) \wedge \star \eta_{(l)}  \hat \varphi_{(l)}\right)
= 
2 \imath g^{\alpha\alpha} \omega_{(k), \alpha}\breve{b}_{(k)} \eta_{(l), \alpha}  \hat{\varphi}_{(l)},
\end{align*}
and the second term is treated analogously. We turn to the third term, and use \eqref{eq : explicit bullet}
\begin{align*}
-\frac{\imath}{2} (\eta_{(k)} \wedge \hat W_{(k)}) 
\bullet \hat \varphi_{(l)}
&=
\frac{\imath}{2} \Gamma^\alpha \Gamma^\beta (\eta_{(k),\alpha} \omega_{(k),\beta} - \eta_{(k),\beta} \omega_{(k),\alpha}) \breve{b}_{(k)} \hat \varphi_{(l)}
\\&=
\imath \Gamma^\alpha \Gamma^\beta \eta_{(k),\alpha} \omega_{(k),\beta} \breve{b}_{(k)} \hat \varphi_{(l)},
\end{align*}
where the second equality follows from \eqref{commrel_Gamma} together with $\eta_{(k),\alpha} \omega_{(k)}^{\alpha} = 0$.
The fourth term is similar, and we obtain
    \begin{align}\label{two_fold_for_code}
 		\hat{N}^{\mathrm{D}}_{(kl)}  
 	&= 2 g^{\alpha\alpha} \omega_{(k), \alpha}\breve{b}_{(k)} \imath \eta_{(l), \alpha}  \hat{\varphi}_{(l)} + 2 g^{\alpha\alpha} \omega_{(l), \alpha}\breve{b}_{(l)}  \imath \eta_{(k), \alpha} \hat{\varphi}_{(k)} 
 	\\\notag&\quad +  \Gamma^\alpha\Gamma^\beta \imath \eta_{(k), \alpha} \omega_{(k), \beta}\breve{b}_{(k)}  \hat{\varphi}_{(l)}  +  \Gamma^\alpha\Gamma^\beta \imath \eta_{(l), \alpha}  \omega_{(l), \beta}\breve{b}_{(l)} \hat{\varphi}_{(k)}. 
	\end{align}
Finally, as $\eta_{(kl)}$ is not light-like, it holds that 
\begin{align}\label{two_fold_phi_for_code}\hat{\varphi}_{(kl)} = \sigma[\Box]^{-1}(y, \eta_{(kl)}) \hat{N}^{\mathrm{D}}_{(kl)}.\end{align}

\subsubsection{Three-fold interactions}

Analogously to the two-fold terms, we expand the three-fold terms
	\begin{align}\label{three_fold_for_code}
\hat N^{\mathrm{D} (1+1+1)}_{(123)} 
&= 
\sum_{\pi \in S_3}
g^{\alpha \alpha} \omega_{(\pi(1)), \alpha}\breve{b}_{(\pi(1))}
\omega_{(\pi(2)), \alpha}\breve{b}_{(\pi(2))}
\hat \varphi_{(\pi(3))}  
\\\notag
\hat N^{\mathrm{D}\wedge}_{(123)} 
&=	 
\sum_{\pi \in S_3}
g^{\alpha \alpha} \omega_{(\pi(1)), \alpha}\breve{b}_{(\pi(1))}
\imath \eta_{(\pi(2)\pi(3)), \alpha} \hat\varphi_{(\pi(2)\pi(3))}
\\\notag
\hat N^{\mathrm{D}\bullet}_{(123)} 
&=  
\frac{1}{2} \sum_{\pi \in S_3}
\Gamma^\alpha \Gamma^\beta 
\imath\eta_{(\pi(1)),\alpha} 
\omega_{(\pi(1)), \beta}\breve{b}_{(\pi(1))}
\hat\varphi_{(\pi(2)\pi(3))}.
	\end{align} 
    
Recall that the notation $I_\gamma$ is defined by \eqref{truncated_ray_trans}, and write $I_{(j)} = I_{\gamma_{(j)}}(\ell)$ where $\gamma_{(j)}$
is as in \eqref{def_gamma_i}. Recall that $\sigma[\varphi_{(j)}]$ satisfies the transport equation \eqref{transport_red} on $\gamma_{(j)}$ and that $\hat \varphi_{(j)}$ is a rescaled version of $\sigma[\varphi_{(j)}]$ at $(y, \eta_{(j)})$. Hence, in view of \eqref{transport_sol},
	\begin{align}\label{hat_phi_for_code}
\hat \varphi_{(j)} = -\frac12 \Gamma^\alpha \omega_{(j),\alpha} \Gamma^\beta \xi_{(j),\beta} \breve b_{(j)} I_{(j)}.
	\end{align} 
Passing to the limit $s \to 0$, the convergence $\xi_{(j)} \to \xi_{(1)}$ entails that $I_{(j)} \to I_{(1)}$ as well.
We define 
	\begin{align*}
I_{(j)}' = \p_s I_{(j)}|_{s = 0}, 
\quad j=2,3.
	\end{align*}

Simplification of $\hat{N}_{(123)}^{\mathrm{D}}$ is in principle straightforward but in practice very tedious. We have verified the following formula using a symbolic computation, see \cite{CLOP2025_code},
	\begin{align*}
-\frac{4}{3 s^2 r} \hat{N}_{(123)}^{\mathrm{D}} 
&=
\frac{1}{r} \breve b^{3} (\Gamma^0 - \Gamma^1) (\Gamma^2  I_{(1)}
+ \Gamma^1 (I_{(2)}' - I_{(3)}'))
\\\notag&\qquad
+ \breve b^{3} I_{(1)} 
- \frac{1}{2} \breve b^{3} (\Gamma^0 + \Gamma^1) \Gamma^2 (I_{(2)}' - I_{(3)}')
+ \mathcal{O}(s) + \mathcal{O}(r),
	\end{align*}
where $\breve b = \varrho_*(b)$, with $b$ as in \eqref{eqn : special sources}.
In particular, as $(\Gamma^0 + \Gamma^1)^2 = 0$, it holds that
	\begin{align}\label{final_interaction_mod}
-\frac{4}{3 s^2 r} (\Gamma^0 + \Gamma^1) (\hat{N}_{(123)}^{\mathrm{D}} - \lim_{r \to 0} \hat{N}_{(123)}^{\mathrm{D}})
&= \breve b^3 (\Gamma^0 + \Gamma^1) I_{(1)} + \mathcal{O}(s) + \mathcal{O}(r).
	\end{align}

\subsection{Retrieval of the Dirac components}\label{subsec : retrieval of Dirac}

Recall that the light ray $\gamma_{(4)}$ from $y \in \mathbb D$ to $z = \gamma_{(4)}(\ell) \in \mho$ is defined by \eqref{def_gamma_4}.
Then for $\sigma[(\varphi_{(123)}, W_{(123)}, \Upsilon_{(123)})] = (\varsigma, w, \upsilon)$ the transport equations \eqref{eq:ymhdyD+}-\eqref{eq:ymhdyH} along $\gamma_{(4)}$ have vanishing initial conditions for $w$ and $\upsilon$, whereas the initial condition for $\varsigma$ is non-vanishing and proportional to $\hat N_{(123)}^{\mathrm{D}}$. 
Hence the Dirac channel component $\hat \varphi_{(123)}$ of
	\begin{align*}
\sigma\left[ \left.\p^3_{\epsilon_{(1)}\epsilon_{(2)}\epsilon_{(3)}} \mathbf{L}_{(\psi, A, \Phi)}(0, \mathcal{J}_1, \mathcal{J}_2, \mathcal{J}_3, \mathcal{F}) \right|_{\epsilon = 0} \right] (z, \eta)
	\end{align*} 
satisfies $\hat \varphi_{(123)} = c \P_{\gamma_{(4)}}^{A,\varrho}(\ell) \hat N_{(123)}^{\mathrm{D}}$,
where $c \ne 0$ does not depend on $\psi$ (but it does depend on $s$).
In particular, we can determine
	\begin{align*}
\hat N_{(123)}^{\mathrm{D}} = c^{-1} \P_{\gamma_{(4)}}^{A,\varrho}(\ell)^{-1} \hat \varphi_{(123)}.
    \end{align*} 

As $\varrho$ is hypercharged, we can choose such $b \in Z(\g)$ that
$\varrho_*(b) = \breve b$ is invertible. 
Then it follows from \eqref{final_interaction_mod} that we can determine
	\begin{align*}
(\Gamma^0 + \Gamma^1) I_{(1)}.
    \end{align*}
As $I_{(1)} = I_{\gamma_{(1)}}(\ell)$, moving $y$ slightly along $\gamma_{(1)}$ 
allows us to recover $v \cdot I_{\gamma_{(1)}}(t)$ for $t$ near $\ell$, where $v = \dot \gamma_{(1)}$. In view of \eqref{truncated_ray_trans},
we can also recover
	\begin{align*}
v \cdot \psi(y) = \P_{\gamma_{(1)}}^{A,\varrho} \p_t (\P_{\gamma_{(1)}}^{A,\varrho})^{-1} (v \cdot I_{\gamma_{(1)}})|_{t=\ell}.
	\end{align*}

We repeat the above construction with the same $y \in \mathbb D$ but with $x \in \mho$ replaced by its small perturbation that is not on the line through $x$ and $y$. 
This amounts to the recovery of $\tilde v \cdot \psi(y)$
for a lightlike vector $\tilde v$ that is linearly independent of $v$. Then we obtain
$(v + \tilde v) \cdot \psi(y)$.
The vector $v + \tilde v$ is not lightlike, as it is the sum of two linearly independent lightlike vectors. Hence it is invertible in the Clifford algebra, and we recover $\psi(y)$.
This concludes the proof of Theorem \ref{thm : StS to fields}.

	\bigskip
 
\noindent {\bf Acknowledgements.} XC was supported by Natural Science Foundation of Shanghai grant 23JC1400501. 
ML was supported by the European Research Council of the European Union (ERC), grant 101097198,
and the Research Council of Finland (RCF), grant 359186.
LO was supported by ERC, grant 101086697, and RCF, grants 347715, 353096 and 359182. GPP was supported by NSF grant DMS-2347868.
Views and opinions expressed are those of the authors only and do not necessarily reflect those of the European Union or the other funding organizations.



\bibliographystyle{abbrv}
\bibliography{main}

\end{document}